\documentclass[12pt,a4paper]{article}
\usepackage{PRIMEarxiv}

\usepackage[utf8]{inputenc}
\usepackage{alphabeta}

\usepackage{todonotes}
\usepackage{geometry} % [top=1in, bottom=1in, left=1in, right=1in]

\usepackage[parfill]{parskip}

\newcommand{\eat}[1]{}

\usepackage{enumitem}
\usepackage{bbm}
\setlist{nolistsep,noitemsep}
\usepackage{multirow}
\usepackage[hidelinks]{hyperref}

\usepackage[greek,english]{babel}
\usepackage{lipsum}
\usepackage{lastpage}
\usepackage{esvect}
\usepackage{fancyhdr}
\usepackage{algorithmic}
\usepackage{subcaption}
\usepackage[numbers]{natbib}
\usepackage{placeins}
\usepackage{booktabs}
\usepackage{xcolor}
\usepackage{mathtools}
\usepackage{wrapfig}
\usepackage{comment}
\usepackage{amssymb}
\usepackage{bm}
\hypersetup{
  colorlinks   = true, %Colours links instead of ugly boxes
  urlcolor     = blue, %Colour for external hyperlinks
  linkcolor    = blue, %Colour of internal links
  citecolor   = blue %Colour of citations
}
\usepackage{float}
\floatstyle{plaintop}
\restylefloat{table}

\usepackage{amsmath,amsfonts}
\usepackage[noend,ruled,vlined]{algorithm2e}
\usepackage{array}

\usepackage{cleveref}

\title{
Bike network planning in limited urban space
}
\author{
Nina Wiedemann$^{1, *}$, Christian Nöbel$^{2}$,  Lukas Ballo$^{3}$, Henry Martin$^{1}$, Martin Raubal$^{1}$
}

\usepackage[some]{background}
\backgroundsetup{%
  scale=1,       %% change accordingly
  angle=0,       %% change accordingly
  opacity=0.7,    %% change accordingly
  color = {darkgray},  %% change accordingly
  contents={\parbox{\textwidth}{\centering %  12cm
 \large This paper has been accepted for publication in \\ Transportation Research Part B: Methodological \\[25cm]}
}
}

\begin{document}

% \tableofcontents
\newpage

\maketitle
\vspace{-2em}
\begin{center}
\small{$^{*}$\textmd{Corresponding author: \texttt{nwiedemann@ethz.ch}}}\\
\small\textsc{$^{1}$ Institute of Cartography and Geoinformation, ETH Zurich, Switzerland}\\
\small\textsc{$^{2}$ Institute for Operations Research, ETH Zurich, Switzerland}\\
\small\textsc{$^{3}$ Institute for Transport Planning and Systems, ETH Zurich, Switzerland}\\
\end{center}
\vspace{4em}

\BgThispage{}

\begin{abstract}

The lack of cycling infrastructure in urban environments hinders the adoption of cycling as a viable mode for commuting, despite the evident benefits of (e-)bikes as sustainable, efficient, and health-promoting transportation modes. Bike network planning is a tedious process, relying on heuristic computational methods that frequently overlook the broader implications of introducing new cycling infrastructure, in particular the necessity to repurpose car lanes. In this work, we call for optimizing the trade-off between bike and car networks, effectively pushing for Pareto optimality. This shift in perspective gives rise to a novel linear programming formulation towards optimal bike network allocation. Our experiments, conducted using both real-world and synthetic data, testify the effectiveness and superiority of this optimization approach compared to heuristic methods. In particular, the framework provides stakeholders with a range of lane reallocation scenarios, illustrating potential bike network enhancements and their implications for car infrastructure. Crucially, our approach is adaptable to various bikeability and car accessibility evaluation criteria, making our tool a highly flexible and scalable resource for urban planning. This paper presents an advanced decision-support framework that can significantly aid urban planners in making informed decisions on cycling infrastructure development.

\end{abstract}


\vspace{2em} \section{Introduction}
\vspace{1em}

The transportation sector plays a pivotal role in combating climate change, accounting for around 20.7\,\% of CO$_2$ emissions worldwide~\citep{owid-ghg-emissions-by-sector}. % Of particular concern is the continuous rise in emissions, which threatens to undermine the progress achieved in other sectors aimed at reducing carbon footprints.
In urban environments, a viable alternative to motorized transport is cycling, promising sustainable traffic in addition to % as a sustainable alternative for urban traffic thus plays an important role.
% In urban environments, cycling presents itself as the most sustainable transport option, 
% additionally 
% Beyond the clear advantage of reducing emissions, cycling also 
substantial health benefits~\citep{oja2011health}. The emergence of e-bikes in recent years has further democratized cycling for the general population. 
Nevertheless, in most cities cycling only accounts for a minority of transport activity so far. 
A crucial factor in the adoption of bicycle commuting is the availability and density of bike networks~\citep{schoner2014missing}. Previous research has provided compelling evidence for the positive impacts of bike lane infrastructure on cycling~\citep{buehler2016bikeway}, encompassing its effects on public health~\citep{mueller2018health}, the importance of physical separation from car lanes~\citep{fraser2011cycling}, and the incorporation of green spaces~\citep{ta2021impact}. These findings collectively underscore the importance of well-designed bike infrastructure in promoting sustainable and healthy urban transportation choices.

Many cities around the world have recently promoted the construction of large-scale cycling infrastructure. Implementation often suffers from practical difficulties and a complex planning process that involves many stakeholders. Computational methods can serve as decision support systems, including topology-based methods aiming to improve the connectivity of bike networks~\citep{natera_orozco_data-driven_nodate}, 
cost-benefit analyses~\citep{paulsen2023societally, szell_growing_2022}, or data-driven planning methods based on data from bike sharing systems~\citep{steinacker_demand-driven_2022, duthie2014optimization, bao2017planning, liu2022urban}, travel surveys~\citep{mauttone2017bicycle} or mobile phones~\citep{olmos2020data}. 
Methods range from simple heuristics, e.g. based on the betweenness centrality of edges in the network~\citep{steinacker_demand-driven_2022, ballo_e-bike_2023}, to linear programming (LP)~\citep{duthie2014optimization, lin2013bikeway} or  mixed integer linear programming (MILP)~\citep{liu2022urban} approaches. 

However, most approaches that propose large scale changes of the existing 
infrastructure do not consider the impact of these changes on other 
transport modes~\citep{gerike2022network}.
As cities have a (mostly) fixed street infrastructure, an improved cycling infrastructure is only attainable by allocating existing streets or lanes from car to bike usage.
Thus, every improvement in bike infrastructure comes at the cost of worsening the car network~\citep{burke2016space}.
As such an essential factor for the practical feasibility and public acceptance of a radical network restructuring in favor of cycling is its impact on the car network.

The contribution of this work is two-fold. 
First we propose a framework to evaluate bike network planning approaches % for essential changes of the car network 
that takes into account both the improvement of the bike infrastructure as well as the car reachability in the modified network. 
This is achieved by comparing planning methods in terms of their \textit{Pareto frontier}, trading off the bikeability against the car travel times.  
% considering the travel times for both modes of transportation over different proposed solutions computed by a framework.
Informally speaking, the solutions are compared based on the question ``How much can one improve the bike network, when allowing the car reachability to deteriorate by no more than a fixed amount''. 
%Different frameworks can then be compared based on how well they can connect the bike network for a fixed worsening of the car network.
Secondly, we propose a novel optimization framework for bike network planning. 
In the spirit of the pareto-optimality-based evaluation framework, our optimization method minimizes both bike and car travel times in a multi-criteria objective function. % when planning the reallocation of lanes.
The approach, based on a flow modeling LP, is flexible both in regards to different trade-offs between car and bike network, as well as for modeling further constraints. 
In particular, the framework can provide stakeholders with a multitude of possible reallocations, showing which improvements of the bike network can be achieved at what cost for the car infrastructure.

% results
Experiments on simulated and real data from several cities demonstrate that our approach yields solutions superior to previous heuristic approaches, in terms of the dominance of the Pareto frontier. In addition, we demonstrate the successful trade-off between computational efficiency and solution quality achieved with the proposed relaxations. The framework can be easily extended to incorporate further constraints such as parking spaces, dedicated bus lanes, pedestrian infrastructure, and unsealed green spaces. Thus, it can serve as a tool to support urban planners in designing future sustainable cities, leveraging cycling with minimal impact on other transportation. 

The remainder of the paper is organized as follows: \Cref{sec:relatedwork} discusses existing literature. In \cref{sec:methods} our LP formulation is explained together with certain relaxations. \Cref{sec:real_data} and \cref{sec:synthetic} present our results on real and synthetic data respectively, and in \cref{sec:case_study} we conduct a case study for the city of Zurich. \Cref{sec:discussion} discusses the findings and concludes the paper. 

\vspace{2em} \section{Related work}\label{sec:relatedwork}
\vspace{1em}

\vspace{1em} \subsection{Related combinatorial problems}\label{sec:relatedwork1}
\vspace{1em}

Bike network planning belongs to a class of optimization problems termed ``Urban network design problem'' (UNDP). \citet{farahani_review_2013} review literature on the UNDP and divide it into different categories, where our problem formulation is most closely related to the ``Discrete Network Design Problem''. \citet{farahani_review_2013} find that these problems are oftentimes solved with a bi-level formulation: the top level optimizes the topological layout, and the lower level determines the effect of the proposed layout on the travel routes of individuals. \citet{gallo_meta-heuristic_2010}, for example, optimize the road directions and signal settings with a bi-level optimization approach. \citet{leblanc1975algorithm} addressed the problem of selecting which links in an urban road network should be improved to minimize total congestion, using a MILP approach. Similarly, \citet{gao2005solution} proposed a classical bi-level formulation for selecting link additions to an existing road network, considering fixed demand from each origin to each destination. 

Another line of work, starting with \citet{foulds1981multi}, framed urban network design as a multi-commodity flow problem~\citep{ford1958suggested}, which provides a powerful framework for traffic modeling. \citet{bevrani2020multi} take this approach to analyze the effect of network changes on traffic flow, following previous work on using multi-commodity flow models as macroscopic traffic simulations~\cite{nilsson2013multi, wright2017node}. Such traffic flow modeling approaches were further leveraged for optimization, e.g. for network routing~\citep{poh2004multicommodity} or signal control~\citep{de2020multi}.

However, bi-level optimization and multi-commodity flow formulations suffer from high computational costs. Due to the discrete decision variables that determine lane directedness or lane existence, network design problems are usually NP-hard. Even the simple problem of optimizing edge directions (via edge flipping), with the goal of minimizing the number of connected components of the graph, was shown to be NP-hard~\citep{kanade_connectivity_2004}. Thus, heuristics and approximations are widely used for network design. \citet{cantarella_multi-criteria_2006} propose a genetic algorithm to determine suitable configurations of traffic lights, parking lots, link layout, and further design choices. \citet{poorzahedy2005application} proposed an ant system for solving network design problems, demonstrating the applicability of bio-inspired algorithms to this domain. \citet{burke2016space} took a planning perspective, suggesting the removal of lanes from motorized metrics based on their potential to disrupt traffic, which involved highly iterative processes to refine the network layout. Additionally, \citet{miandoabchi2012bi} employed a genetic programming approach to optimize bus networks, further illustrating the diverse methodologies applied in network design problems. These approaches are promising starting points but need to be adapted to account for the unique challenges of bike network design.

\vspace{1em} \subsection{Heuristic and optimization approaches to bike lane allocation}
\vspace{1em}

Due to the complexity of optimization methods for network design, most previous methods for bike network planning are based on heuristics or iterative approaches for \textit{growing} networks, in contrast to attempts of finding globally optimal networks. We provide an overview in \autoref{tab:literature}, distinguishing previous works by method (i.e. optimization-based vs heuristics), by whether they are data-driven (e.g. incorporating bike sharing data) and by the objective functions used. On the one hand, there are optimization approaches that usually operate in a limited setting, such as allocating bike lanes between bike sharing stations, via an LP~\citep{duthie2014optimization}, IP~\citep{lin2013bikeway},  MILP~\citep{liu2022urban} or genetic programming~\citep{martinez2016cycling} formulations. \citet{mauttone2017bicycle} optimize the travel times along origin-destination (OD) paths via a multi-commodity network flow problem. \citet{bao2017planning} show that their formulation of bike-sharing-trajectory-coverage is NP-hard and propose a heuristic. Due to the complexity of optimization methods, the size of the inputs is oftentimes limited in these studies. Nevertheless, they can be applied to real-world cases after simplifications of the road network; for example, \citet{lin2013bikeway} optimize an area in Taipei City with 75 nodes and 115 candidate links.

On the other hand, there are purely heuristic approaches that are, in turn, applicable to large-scale street networks. Recently, \citet{steinacker_demand-driven_2022} propose a demand-driven heuristic algorithm, starting from a full bike network and generating a sequence of networks by iteratively removing edges with the lowest betweenness centrality. \citet{ballo2023} integrate this approach into their pipeline for bike network planning, but essentially revert the algorithm by starting from a complete \textit{car} network and iteratively reallocating lanes to the bike network, starting from the ones with minimal betweenness centrality in the car network. Other heuristic methods target new objectives, such as improving the connectedness of existing bike networks with new bike lanes~\citep{natera_orozco_data-driven_nodate}, or maximizing \textit{trajectory} coverage~\citep{bao2017planning}. \citet{paulsen2023societally} take a fundamentally different approach by optimizing the cost-benefit trade-off for the sake of long-term planning. 

An important aspect of bike network planning is the demand modeling, since cycling route choices differ significantly from motorized traffic~\citep{meister2023route}. Recent approaches  tend to be data-driven and tap into new data sources, such as \citet{olmos2020data} who leverage mobile-phone-data to assess travel demand, or \citet{lin2013bikeway} who incorporate traffic accident reports to assess cycling safety. In addition to census data and bike sharing data, these data sources are highly relevant to paint a comprehensive picture of cycling behavior.
%
% many focused on bike sharing stations --> limited scope

A common limitation of previous work is its neglect of the impact on motorized traffic. \citet{gerike2022network} review bike network design approaches and conclude that ``the main challenge that is hardly addressed in any of the identified references is the coordination of cycle network development with the 
other transport modes and street functions'' (\citet{gerike2022network}, page 78). % Indeed, we argue that the trade-off between bike and car network has not been addressed sufficiently in previous work.
There are notable exceptions, such as \citet{mesbah2012bilevel} who minimize a weighted sum of the travel times for bicyles and motorized traffic with a bi-level formulation based on a genetic algorithm. However, a proper formulation of the problem, as well as a flexible optimization framework that scales to real-world instances is still lacking. 

Our framework for bike network planning closes this gap by proposing a multi-criteria optimization approach. At the same time, we combine the advantages of previous methods, as the optimal planning method should 1) be optimization-based to provide a global perspective on the network planning task and avoid local minima, 2) be efficient and scale to real-world street networks instead of focusing on bike sharing station networks, 3) be data-driven, i.e. be able to incorporate various travel demand data such as OD-data, census data or bike sharing trajectories.

\begin{table}[h!]
\centering
\caption{Summary of existing literature on bike lane allocation}
\label{tab:literature}
\resizebox{\textwidth}{!}{
\begin{tabular}{|l|l|l|l|}
\hline
\textbf{Paper} & \textbf{Method} & \textbf{Data} & \textbf{Objective} \\ \hline
\citet{paulsen2023societally} & IP & Cycling network Copenhagen & Cost-benefit trade-off \\ \hline
\citet{mesbah2012bilevel} & Bi-level, genetic algorithm & Synthetic data &
\begin{tabular}[t]{@{}l@{}} Minimize travel times for bicycles \\ and motorized traffic
\end{tabular}
\\ \hline
\citet{duthie2014optimization} & LP & Bike sharing & \begin{tabular}[t]{@{}l@{}} Minimize costs for covering \\ demand (bound bike travel times)  \end{tabular} 
\\ \hline
\citet{lin2013bikeway} & Multi-objective IP & Traffic accident data & 
\begin{tabular}[t]{@{}l@{}} Maximize cycling safety \& comfort, maximize \\ coverage, minimize impact on other transport \end{tabular} 
% case study with 75 nodes and 115 candidate links (but still say that it covers 5km radius area
% impact on other transport is minimized by adding an objective that penalizes reduction of parking space and reduction of space for other driving modes
\\ \hline
\citet{mauttone2017bicycle} & Multi-commodity  flow MIP & OD data (surveys) & Cost-benefit trade-off (route coverage) \\ \hline
\citet{liu2022urban} & MILP & Bike sharing data & Maximize coverage of route choices \\ \hline
\citet{martinez2016cycling} & Genetic programming & Bike rental data & Shortest paths between stations \\ \hline
\citet{steinacker_demand-driven_2022} & Heuristic (betweenness) & Bike sharing data + OSM & Minimize bike travel times \\ \hline
\citet{ballo2023} & Heuristic (betweenness) & Census data + OSM & Minimize impact on car network \\ \hline
\citet{bao2017planning} & Heuristic & Bike sharing trajectories &
\begin{tabular}[t]{@{}l@{}} Maximize trajectory coverage \\ and connectivity at minimal costs
\end{tabular}
\\ \hline
\citet{natera_orozco_data-driven_nodate} & Heuristic (add links greedily) & OSM & Improve network connectedness \\ \hline
\citet{olmos2020data} & Percolation theory & Mobile-phone data + census & Maximize travel demand coverage \\ \hline
\end{tabular}
}
\end{table}

\vspace{2em} \section{Methods}\label{sec:methods}

\vspace{1em} \subsection{Problem definition}\label{sec:problemdefinition}
\vspace{1em}
We propose a multi-modal view on bike network planning that considers the impact of new bike lanes on the car network. It is based on the following assumptions: 1) The input is a given street network of an urban area. It is not possible to build entirely new infrastructure, but the type and division of existing road space can be changed, including lane directions. 2) The allocation of bike lanes inevitably involves the reduction of street space available to other modes. 

This planning problem, here termed the ``bike network allocation problem'' (BNAP), can be modeled as a graph division problem. 
The initial graph is a simplified version of the existing street network, corresponding to an undirected graph $G = (V, E)$, where the nodes $V$ are intersections and the edges $E$ are streets between two intersections. Geographic properties of the network are expressed in attributes of the street edges $e$, specifically the length of the street $d(e)$, its gradient $\delta(e)$, the speed limit $\theta(e)$ and its capacity $\Lambda(e)$. The capacity can be set to the street width or the number of lanes. 

Solving the BNAP involves dividing $G$ into two graphs: the bike lane network $G_b = (V, E_b)$ and the car lane network $G_c = (V, E_c)$. Both are directed multi-graphs, since their edges now represent \textit{lanes}. The car network $G_c$ must be strongly connected; i.e., every node must be reachable from any other node, since disconnected subgraphs are unrealistic in an urban environment. The design of $G_b$ and $G_c$ is mainly constrained by the street capacities, ensuring that the car and bike lanes fit into the existing infrastructure:
\begin{align}\label{eq:bnap_constraint}
\forall e=(u,v) \in E: \lambda^c_{(u,v)} + \lambda^c_{(v, u)} + 0.5 (\lambda^b_{(u, v)} +  \lambda^b_{(v,u)}) \leq \Lambda(e) \
\end{align}
 Here, $\lambda^c_{(u,v)}$ denotes the capacity of a directed edge $e=(u,v)$ allocated for cars and % $\Lambda(e), e\in E_c$ and
$\lambda^b_{(u,v)}$ the capacity allocated for cycling. If the edge $(u,v)$ is not part of the car network, we define $\lambda^c_{(u,v)}=0$ and accordingly for the bike network. \autoref{eq:bnap_constraint} defines a division of the total capacity of the undirected street into directed car and bike lanes. The bike capacities are multiplied by 0.5 to express the lower space necessary for bikes. We follow common guidelines for bicycle infrastructure that recommend a bike lane width of 1.5 m~\citep{parkin2018designing, yan2018recommended}, corresponding to about half of a car lane.

\vspace{1em} \subsection{An evaluation framework based on Pareto optimality}\label{sec:bnap}
\vspace{1em}

While a plethora of metrics has been proposed for evaluating ``bikeability''~\citep{weikl_data-driven_2023, grisiute2023ontology}, they largely ignore the impact of proper bike lanes on motorized travel, with the exception of \citet{burke2016space} who propose the Network Robustness Index to measure the effect of wider bike lanes on traffic. For simultaneously evaluating the goodness of the bike and car network, we leverage the concept of Pareto optimality, referring to solutions where improving one aspect (e.g., the bike travel time) inevitably results in the worsening of another aspect (e.g., the car travel time) and the other way around. Pareto optimality has only been used in related work to reflect a multi-criteria evaluation of bikeability with respect to distance, safety, comfort, etc~\citep{reggiani2022understanding, gholamialam2019modeling}. Instead, we compute the Pareto frontier between car and bike travel times, since travel time arguably is a main indicator of the goodness of a street network. 

\begin{table}[ht]
    \centering
    \begin{tabular}{rr}
    Variable & Description \\
    \toprule
    \textbf{Optimized variables} & \\
      $f_{s,t,e}^c$   &  Car flow on car lane along edge $e$ for the path from $s$ to $t$ \\
      $f_{s,t,e}^b$   &  Bike flow on bike lane along edge $e$ for the path from $s$ to $t$  \\      
      $f_{s,t,e}^{\beta}$   &  Bike flow on car lane (shared lane) along edge $e$ for the path from $s$ to $t$  \\
      $\lambda^c_e$ & Capacity allocated for cars on edge $e$ \\       
      $\lambda^b_e$ & Capacity allocated for bikes on edge $e$ \\[2ex]         
      \textbf{Inputs} & \\
      $G = (V, E)$ & Street graph with one edge per street\\
      $G' = (V, E')$ & Auxiliary street graph with a pair of reciprocal directed edges per street\\
      $n$ & Number of nodes\\
      $m$ & Number of edges\\
      $\Omega$ & Set of considered origin-destination pairs, given as pairs of node $(u,v)\in V^2$ \\
      $t^b(e)$ & Travel time by bike along edge $e$ \\
      $t^c(e)$ & Travel time by car along edge $e$ \\
      $t^{\beta}(e)$ & Perceived travel time for cycling on a car lane along edge $e$ \\
      $\gamma$ & Weighting (desired importance) of the car travel time \\
      $\omega(s, t)$ & Travel demand (importance) of the route from $s$ to $t$ \\
      $\lambda_e$ & Overall given capacity of edge $e$ \\ 
      $\phi(s, t)$ & Desired flow from $s$ to $t$
    \end{tabular}
    \caption{List of variables and explanation}
    \label{tab:var_description}
\end{table}

Optimally, the travel times should be computed in agent-based mesoscopic simulations to reflect traffic flow and the agents' mode choices. However, simulations are computationally expensive and impede efficient evaluation of a large number of potential networks. Thus, we compute network-based travel times, i.e., weighted shortest paths, but build on previous work from \citet{steinacker_demand-driven_2022} by taking a demand-driven view on street network goodness. To reflect real-world demand, we incorporate an origin-destination (OD) matrix. Let $\Omega$ be a set of OD pairs where the origin and destination are nodes in the graph, $\Omega = \{(u_1,v_1), (u_2, v_2), \dots \}$. $\Omega$ is derived from travel surveys, GPS trajectories, or bike sharing data. Let $t^c(e)$ and $t^b(e)$ be the travel times by car or bike along edge $e$, and let $P^c(u, v)$ and $P^b(u, v)$ be the edge set of a weighted shortest path from $u$ to $v$, based on the edge weights $t^c$ and $t^b$ respectively. The Pareto frontier is computed to trade-off the goodness of the car network, $\mathcal{T}(G^c)$, with the goodness of the bike network $\mathcal{T}(G^b)$, defined as
\begin{align}\label{eq:pareto_objective}
    \mathcal{T}(G^c) = \sum_{(u,v)\in \Omega} \Big( \sum_{e\in P^c(u, v)} t^c(e) \Big)\ \ \ \ \   \ \ \ \   \ \ \ \    \mathcal{T}(G^b) = \sum_{(u,v)\in \Omega} \Big( \sum_{e\in P^b(u, v)} t^b(e) \Big)
\end{align}

The edge-wise travel times $t^b(e)$ and $t^c(e)$ are set based on the lane's length $d(e)$ in km, its gradient $\delta(e)$ in \% and its speed limit $\theta(e)$ in km/h, which are available from Open Street Map (OSM) data. For the car travel times we simply set $t^c(e) =\frac{d(e)}{\theta(e)}$ since the gradient does not have a strong impact on the car speed in urban areas. 
Bike speed on the other hand is highly impacted by the gradient. This effect is estimated by \citet{parkin2010design}. They determine the base speed of cyclist to be 21.6km/h (excluding acceleration and breaking) and they further derive an increase of 0.86km/h per negative percent gradient (downhill acceleration) and an increase of 1.44km/h per positive percent gradient (uphill). For example, a lane with 3\% negative gradient results in an estimated bike speed of % $v^b(e) = 16.35 + 0.86 \cdot 3 = 18.93$.
$v^b(e) = 21.6 + 0.86 \cdot 3 = 24.18$ km/h. Together, $t^b(e)$ is set as follows:
\begin{align}
    t^b(e) = \frac{d(e)}{v^b(e)} = \begin{cases}
        \frac{d(e)}{max\{1,\ \ 21.6  - 1.44\cdot \delta(e)\}} \ & \text{if}\ \delta(e) > 0 \text{  (uphill)} \\[1.5ex]
        \frac{d(e)}{21.6 - 0.86\cdot \delta(e)} \ & \text{else} \text{  (downhill or flat)}.
        \end{cases}
\end{align}

However, a crucial component of cycling is safety, as a lack of dedicated infrastructure prevents people from cycling. This can be expressed in terms of a \textit{perceived} bike travel time that amends the physical travel time with a psychological component. The perceived bike travel time is computed from the actual travel time, $t^{b}(e)$, by penalizing the discomfort of cycling on car lanes. The penalty is zero for edges equipped with dedicated bike lanes. To set the penalty for roads without proper bike lanes, we follow the study by \citet{meister2023route}, who recently studied route choices of cyclists in Zurich and found a ``value of distance'' of $-0.66$ for bike lanes and $-0.36$ for bike paths. Taking the average of both, we assume that the perceived distance is approximately halved on proper bike lanes, and we thus set perceived travel time for cycling on car lanes to $t^{\beta}(e) = 2 \cdot t^b(e)$.\\ Finally, it is worth noting that our framework could easily be adapted to regard e-bikes by reducing the effect of the gradient on the travel time.

\vspace{1em} \subsection{A linear programming approach towards solving the BNAP}
\vspace{1em}

The goal of minimizing both bike and car travel times gives rise to an integer programming (IP) formulation with a multi-criteria objective function. To express travel times in an IP, we first revisit a flow formulation of the all-pairs shortest path problem:
\begin{align}\label{eq:allpairs1}
% objective
\min \sum_{s,t \in V} \sum_{e\in E}  f_{s, t,e} \cdot t(e)
\end{align}
subject to
\begin{align}\label{eq:allpairs2}
% flow greater zero
% \forall v\in V,\ \ \forall s,t\in V, \forall e\in E: \ \ f_{s,t,e} \geq 0 \\
% flow integer
\forall s,t\in V, \forall e\in E: \ \ f_{s,t,e} \in \mathbb{Z}_{\geq 0} 
% flow constraint car
\end{align}
\begin{align}\label{eq:allpairs3}
\forall v\in V,\ \ \forall s,t\in V: \ \ \sum_{e\in \delta^+(v)} f_{s, t,e} - \sum_{e\in \delta^-(v)} f_{s, t,e} = 
\begin{cases}
-1\ \ \ &\text{if } v=t \\
1\ \ &\text{if } v=s\\
0\ \ & \text{else,}
\end{cases}
\end{align}
where $t\colon E \rightarrow \mathbb{R}_{\geq 0}$ encodes the travel times along the edges, and $f_{s,t,e}$ is the flow allocated on edge $e$ for the path from $s$ to $t$. Additionally $\delta^+(v)$ denotes the set of outgoing edges of node $v$ and $\delta^-(v)$ its incoming edges. The flow constraints, together with the integer constraint $f_{s,t,e} \in \mathbb{Z}$, guarantee that there is a flow of value 1 between every $(s, t)$-pair, corresponding to a path. The objective computes the total travel time along each path, leading to the all-pairs shortest path in $G$. 

We adapt the objective and constraints to model the BNAP. First, the undirected street graph $G$ is converted into a directed graph $G'$ by replacing every edge $e=(u,v)\in E$ by a pair of reciprocal directed edges, $\overrightarrow{e}=(u,v)$ and $\overleftarrow{e} = (v,u)$. The construction of $G'$ allows to optimize the capacity for bike and car lanes in both directions. 
The edge properties of $\overrightarrow{e}$ and $\overleftarrow{e}$ are inherited from $e$.
To make this specific we let the length of the street and its speed limit agree with the original version of the edge, i.e.\ $d(\overrightarrow{e}) = d(\overleftarrow{e}) = d(e)$ and $\theta(\overrightarrow{e}) = \theta(\overleftarrow{e}) = \theta(e)$, respectively.
The gradient attribute is set based on the definition of $\delta(e)$, with $\delta(\overrightarrow{e}) = -\delta(\overleftarrow{e})$.
It is worth noting that we do not assign capacities to the new arcs, as they both compete for the same fixed capacity $\Lambda(e)$.
Instead we will introduce variables for the capacities of $\overrightarrow{e}$ and $\overleftarrow{e}$, together with capacity constraints.

In the following, the adaptation of \autoref{eq:allpairs1}-\autoref{eq:allpairs3} to the BNAP is explained in detail. An overview of all variables is given in \autoref{tab:var_description}.

\subsubsection{Introducing bike, car and shared flow}

In contrast to the general shortest path formulation, we distinguish between the car flow $f^c$ and the bike flow $f^b$ along every edge, and frame the objective as a multi-criteria optimization problem. A weight $\gamma$ can be set by the user to specify the desired importance of the car travel times relative to the bike travel times, resulting in the following preliminary objective function
\begin{align}\label{eq:prelim_time_1}
\min \sum_{e\in E'} \sum_{(s,t)\in \Omega} f_{s, t,e}^b t^b(e) + \gamma \cdot f_{s, t,e}^c t^c(e)
\end{align}
Varying $\gamma$ yields Pareto-optimal solutions; i.e., scenarios where improving the bike network increases the car travel times and the other way around. 
Importantly, we compute the travel times over the s-t-pairs in an OD-matrix~$\Omega$, as motivated in \autoref{sec:bnap}, instead of considering all node pairs. Additionally considering an OD-matrix reduces the runtime significantly. Specifically, with $n = |V|$ and $m=|E'|$, the all-pairs-shortest path formulation requires to optimize $n^2m$ bike and car flow variables, since $n^2$ s-t-paths are considered, each with one variable per edge. This reduces to $|\Omega| m$ variables with the demand-driven approach.

The car flow constraints are set as in \autoref{eq:allpairs3} of the all-pairs shortest-path formulation:
\begin{align}\label{eq:flowcar}
% flow constraint car
\forall v\in V,\ \ \forall (s,t)\in \Omega: \ \ \sum_{e\in \delta^+(v)} f_{s, t,e}^c - \sum_{e\in \delta^-(v)} f_{s, t,e}^c = 
\begin{cases}
-1\ \ &\text{if } v=t \\
1\ & \text{if } v=s\\
0\ \ \ \ &\text{else.}
\end{cases}
\end{align}
However, this constraint is unsuitable for bike flow, since simultaneously guaranteeing bike and car flow along each path in $\Omega$ is oftentimes infeasible in a real-world street network. 
Thus, we introduce the concept of \textit{shared flow}, denoted $f^{\beta}$, representing that cyclists can also use car lanes. \autoref{eq:flowbike} expresses that a bike paths between all pairs of nodes $s,t$ is only required with a combination of $f^b$ (bike-on-bike-lane flow) and shared flow $f^{\beta}$ (bike-on-car-lane flow). This modification ensures that the problem is never infeasible.
\begin {align}\label{eq:flowbike}
% flow constraint bike
\forall v\in V,\ \ \forall (s,t)\in \Omega:\ \ \sum_{e\in \delta^+(v)} (f_{s, t,e}^b + f_{s, t,e}^{\beta}) - \sum_{e\in \delta^-(v)} (f_{s, t,e}^b + f_{s, t,e}^{\beta}) = 
\begin{cases}
-1 & \text{if } v=t \\
1 &\text{if } v=s\\
0 &\text{else.}
\end{cases}
\end{align}

\subsubsection{Constraining the space for bike and car lanes}

Since the space on urban streets is limited, the goal of our approach is to decide which space to allocate to car and bike travel respectively. We model the space limitation with bike and car capacities, denoted  % The bike and car flow along one edge is limited by the respective capacity, whereas shared flow is 
$\lambda_e^b$ and $\lambda_e^c$. The bike and car flow are constrained by the respective capacities:
\begin{align}\label{eq:carcap}
% capacity constraints
\forall (s,t)\in \Omega, \forall e\in E':\ \ f_{s,t,e}^c \leq \lambda^c_e
\end{align}
\begin{align}\label{eq:bikecap}
\forall (s,t)\in \Omega, \forall e\in E':\ \ f_{s,t,e}^b \leq \lambda^b_e
\end{align}

In turn, $\lambda_e^b$ or $\lambda_e^c$ are bounded by the overall street capacity $\lambda_e$ which is given by the number of available lanes or the street width. In alignment with the BNAP (\autoref{eq:bnap_constraint}), it is required that
\begin{align}\label{eq:overallcap}
\forall e=(u,v) \in E': \lambda^c_{(u,v)} + \lambda^c_{(v, u)} + 0.5 (\lambda^b_{(u, v)} +  \lambda^b_{(v,u)}) \leq \Lambda(e) \
\end{align}

Notably, the shared flow $f^{\beta}$ is not constrained, explicitly allowing cycling on any street and in both directions. This is justified since cycling along car lanes is generally permitted, oftentimes including the wrong way along one-way streets. Instead of constraining this undesired bike travel on car lanes, we penalize it by means of a higher perceived travel time $t^{\beta}(e)$. 

The final objective function and full linear program thus becomes:
\begin{gather}\label{eq:final_objective}
 \min \sum_{(u,v)\in E'} \sum_{(s,t)\in \Omega} f_{s, t,e}^b t^b(e) + f_{s, t,e}^{\beta} t^{\beta}(e) + \gamma \cdot f_{s, t,e}^c t^c(e) \\
 \text{subject to}\\
% flow constraint car
\forall v\in V,\ \ \forall (s,t)\in \Omega: \ \ \sum_{e\in \delta^+(v)} f_{s, t,e}^c - \sum_{e\in \delta^-(v)} f_{s, t,e}^c = 
\begin{cases}
-1\ \ &\text{if } v=t \\
1\ & \text{if } v=s\\
0\ \ \ \ &\text{else.}
\end{cases}\\
% flow constraint bikes
\forall v\in V,\ \ \forall (s,t)\in \Omega:\ \ \sum_{e\in \delta^+(v)} (f_{s, t,e}^b + f_{s, t,e}^{\beta}) - \sum_{e\in \delta^-(v)} (f_{s, t,e}^b + f_{s, t,e}^{\beta}) = 
\begin{cases}
-1 & \text{if } v=t \\
1 &\text{if } v=s\\
0 &\text{else.}
\end{cases}\\
% capacity constraints
\forall (s,t)\in \Omega, \forall e\in E':\ \ f_{s,t,e}^c \leq \lambda^c_e\\
\forall (s,t)\in \Omega, \forall e\in E':\ \ f_{s,t,e}^b \leq \lambda^b_e\\
\forall e=(u,v) \in E': \lambda^c_{(u,v)} + \lambda^c_{(v, u)} + 0.5 (\lambda^b_{(u, v)} +  \lambda^b_{(v,u)}) \leq \Lambda(e)
\end{gather}
\begin{gather}\label{eq:end_final}
 f_{s,t, e}^b \geq 0,\ \ \  f_{s,t, e}^c \geq 0, \ \ \  f_{s,t, e}^{\beta} \geq 0  \ \ \  \lambda_e^{c} \geq 0,  \ \ \  \lambda_e^{b} \geq 0
\end{gather}

Solving the problem yields the optimal bike and car capacities $\lambda^c, \lambda^b$ per street and direction. The capacities indicate what fraction of the street width should be allocated to bikes and cars, respectively. If solved as an \textit{integer problem}, the solution can be interpreted as the number of lanes to allocate per transport mode.

\subsubsection{Placing bidirectional bike lanes}
Finally, we add a constraint ensuring that every bike lane is bidirectional:
\begin{align}
\forall e\in E':\ \ \lambda^b_{\overrightarrow{e}} = \lambda^b_{\overleftarrow{e}}
\end{align}
This constraint is motivated by the integer-valued capacities of streets in our dataset, where single-bike-lanes are undesirable due to the infeasibility of half car lanes.\footnote{The capacity $\lambda_e$ is given in terms of the number of lanes and is, therefore, integer-valued (even existing bike lanes are usually bidirectional and together have the capacity of one car lane). The premise of the BNAP is that all available infrastructure is allocated either as bike or car lanes. When placing a single one-directional bike lane on a street with integer-valued capacity, it is only plausible to add another bike lane in the opposite direction, since the remaining half-lane cannot fit a car lane. Therefore, it is more efficient to directly optimize the allocation of bidirectional bike lanes.} However, the constraint is entirely optional and can be easily dropped to extend the formulation to street-width-based bike lane allocation.

\subsubsection{Guaranteeing connectivity via auxiliary paths}\label{sec:methods_od_matrix}

Problematically, omitting the requirement of all-pairs shortest paths, i.e. using an OD-matrix, also removes the requirement of strong connectivity, which is necessary for the car network. 
Therefore, we design $\Omega$ to contain auxiliary $(s,t)$-pairs, which are constructed in chained form to ensure connectivity with the lowest number of variables. In detail, let $v_1, \dots, v_n$ be a random shuffling of all nodes in $V$. We amend the-real-world OD-matrix $\Omega_{demand}$ by the auxiliary pairs $\Omega_{aux} = \{(v_1, v_2), (v_2, v_3), \dots, (v_{n-1}, v_n), (v_n, v_1)$\}, setting $\Omega = \Omega_{demand} \cup \Omega_{aux}$. This ensures strong connectivity, since 1) the flow constraints require a path between every OD-pair in $\Omega$, and 2) any two nodes are connected via the chain of auxiliary nodes. Since $|\Omega_{aux}| = n$, the computational complexity remains linear in $n$, in contrast to the all-pairs shortest path formulation with $n^2$ paths. 

However, the travel times along the auxiliary paths should not affect the objective value, since the paths are not constructed from real-world demand. Thus, we introduce a path-weighting $\omega$ in the objective function, where $\omega^b(s, t)$ and $\omega^c(s, t)$ can be viewed as the real-world importance of the $s-t-$path for cyclists and car drivers respectively.
\begin{align}\label{eq:objective_final_weighted}
\min \sum_{e\in E'} \sum_{(s,t)\in \Omega} \omega^b(s, t) \cdot \big(f_{s, t,e}^b t^b(e) + f_{s, t,e}^{\beta} t^{\beta}(e) \big) +
\gamma \cdot \omega^c(s,t) \cdot f_{s, t,e}^c t^c(e) 
\end{align}

To ignore the auxiliary pairs in the objective, we set $\omega^c(s,t)=\omega^b(s,t)=0\ \ \forall (s,t)\in\Omega_{aux}$. All other s-t-pairs are weighted by 1 in our experiments; however, $\omega$ could be set by the real-world travel frequency between points. It is also worth noting that our formulation could easily be adapted to account for different real-world demand between cyclists and car drives, by distinguishing $\Omega^b$ and $\Omega^c$. For example, $\Omega^b$ could be set based on bike sharing data and $\Omega^c$ based on travel surveys among car drivers.

Another possibility to reflect the real-world importance of certain paths is to require a flow larger than 1 from $s$ to $t$, say $\phi(s,t)$. The flow constraint becomes
\begin{align}\label{eq:phi_traffic_flow}
% flow constraint car
\forall v\in V,\ \ \forall (s,t)\in \Omega: \ \ \sum_{e\in \delta^+(v)} f_{s, t,e}^c - \sum_{e\in \delta^-(v)} f_{s, t,e}^c = 
\begin{cases}
-\phi(s, t) & \text{if } v=t \\
\phi(s, t) & \text{if } v=s\\
0& \text{else}
\end{cases}
\end{align}
In other words, $\phi$ determines the number of paths necessary between each s-t-pair, and can be set based on the travel demand. However, the edge capacity $\lambda_e$ must be adjusted accordingly, and it is unclear how much the capacity must be increased to avoid infeasibility. Thus, we set $\phi(s,t) = 1\ \ \forall s,t\in V$ in our experiments and leave the possibility of modifying the required flow for future work.

\subsubsection{Spatial considerations for runtime improvements}\label{sec:method_kpaths}

\begin{figure}

\centering

    \includegraphics[width=0.6\textwidth]{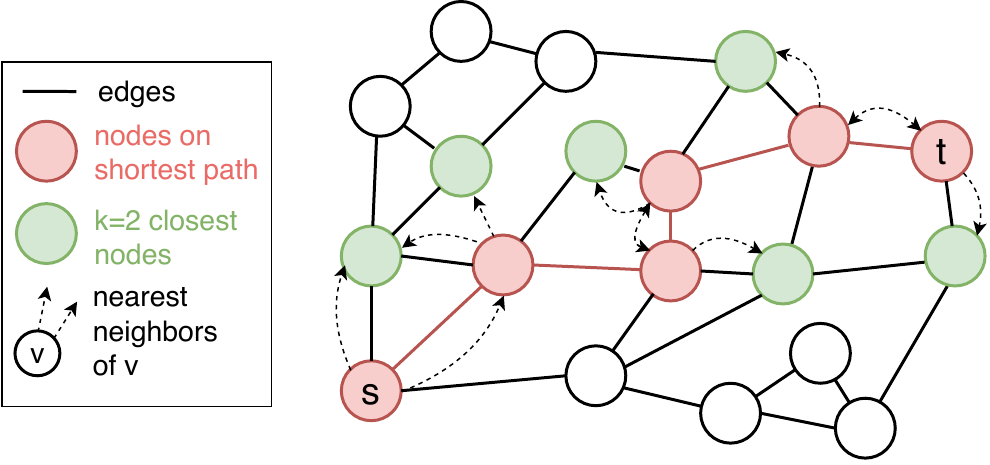}
    \caption{Explanation of spatial relaxation, where flow is only allowed on the subgraph surrounding the shortest path. For each node on the shortest path between s and t, its $\eta$ closest nodes are included. Here closeness is defined in terms of the Euclidean distance in geographic space.}
    \label{fig:kshortest_explained}
\end{figure}

The total number of flow variables in the final formulation is $3|\Omega|\cdot m$, which is still large for real-world instances. To improve runtime efficiency, we propose a relaxation reducing the factor of $m$. In the formulation provided in \autoref{eq:final_objective} -- \autoref{eq:end_final}, one flow variable per $s$-$t$-path and per edge is defined to allow the algorithm to find arbitrary paths through the whole network. However, it is expected that the flow is placed preferably along or close to the shortest path between $s$ and $t$ in order to minimize the overall objective. Based on this assumption, we suggest limiting the set of considered edges per OD-path for large graph instances. Let $\mathcal{E}(s,t)$ denote the set of edges considered by the IP for determining the $s$-$t$ path; i.e., the edges where flow variables are defined.
By default, $\mathcal{E}(s,t)=E'$ for all $(s,t)$-pairs.
We propose to alternatively set $\mathcal{E}(s,t)$ based on the spatial proximity of nodes to the shortest path. Let $P(s,t)=\{v_1, \dots, v_p\}$ be the vertices of a shortest path from $s$ to $t$. Let $\mathcal{N}(v, \eta)$ be the spatial $\eta$-neighborhood of a node $v$, i.e., the set of $\eta$ vertices that is closest to $v$ in geographic space. We define $V_{\eta}(s,t)$ as the set of nodes within the $\eta$-neighborhood of the shortest-path-nodes (see \autoref{fig:kshortest_explained} green), formally $V_{\eta}(s,t) = \bigcup_{v\in P(s,t)} \mathcal{N}(v, \eta)$. We reduce the input graph $G'$ to a subgraph $H$, where $H$ is the induced subgraph based on the node set $V_{\eta}$: $H(s,t) = G'[V{\eta}(s,t)]$.
% We construct a subgraph $H(s,t)=(V_H(s,t), E_H(s,t))$ from $G'$ where , such that the nodes of $H(s,t)$ comprise all nodes. $E_H$ is the set of edges that remains in the induced subgraph. 
When considering only the edges of $H(s,t)$ in the algorithm, i.e. $\mathcal{E}(s,t) = E(H(s,t))$, instead of setting $\mathcal{E}(s,t)=E'$, the number of flow variables per s-t-path is upper bounded by $|P(s,t)| \cdot \eta$.
Usually this bound is quite loose and the actual number of flow variables is much smaller due to overlapping $\eta$-neighborhoods. A visual explanation is provided in \autoref{fig:kshortest_explained}.

\vspace{1em} \subsection{Linear relaxation}
\vspace{1em}

The solution of the provided problem formulation (\autoref{eq:final_objective} -- \autoref{eq:end_final}) can only be translated into bike and car lane-allocations if the problem is solved as an IP. However, the number of variables to optimize remains large even with the proposed relaxations\footnote{Specifically, the total number of variables is $ 3 \cdot m \cdot |\Omega| $ flow variables in addition to $2 \cdot m$ capacity variables ($\lambda^c$ and $\lambda^b$) }, rendering integer programming computationally prohibitive for real-world instances. We propose an iterative process of rounding the capacities and re-computing the optimal solution (\Cref{alg:rounding}).

\begin{algorithm}[htb]
   \caption{Post-processing scheme to round the LP solution}
   \label{alg:rounding}
    \begin{algorithmic}
       \STATE \textbf{Input}: Street network $G$ with edge attributes $d(e), \delta(e), \theta(e)$
       \STATE \textbf{Input}: Set of fixed capacities  $\Lambda$      \STATE \textbf{Input}: $\Omega, \omega, \gamma$, $k$
       \STATE $i = 0$
       \REPEAT
       \STATE $\lambda^* = LP(G, \Lambda, \Omega, \omega, \gamma$)
       \STATE Sort $\lambda^*_b$ % only bikes
       \tcp*{Round edges with largest bike capacity}
       \FOR{k iterations}
            \STATE $\hat{e} = $ edge with largest $\lambda^*_b$ that is not in $\Lambda$ yet
            \IF{$G$ remains strongly connected}
                \STATE $\Lambda = \Lambda \cup \{(\hat{e}, b, 1)\}$ \tcp*{Fix as bike lane}
            \ELSE
                \STATE $\Lambda = \Lambda \cup \{(\hat{e}, c, 1)\}$ \tcp*{Fix as car lane}
            \ENDIF
       \ENDFOR
       \STATE $i = i+1$
       \UNTIL{$|\Lambda| =  2m$} \tcp*{Until all edges are fixed as bike or car lanes}
    \end{algorithmic}
\end{algorithm}

Let $\Lambda$ be a set of the indices and values of all \textit{fixed} capacities, $\Lambda = \{(e, i, \lambda_e^i)\ |\ i\in\{b, c\},\ e\in E\}$. Initially, $\Lambda$ is empty ($\Lambda = \varnothing$), or $\Lambda$ corresponds to a set of lanes that are fixed due to real-world constraints such as compulsory car lanes that are used by bus services. In each iteration, the LP is solved subject to the fixed bike capacities $\Lambda$, yielding the optimal capacities $\lambda^*$, where $\forall (i, e, \lambda_e^i)\in \Lambda:\ \ \lambda_e^i=\lambda_e^{*,i}$ (i.e., the fixed capacities remain unchanged). The algorithm then rounds up the $k$ largest \textit{bike} capacity values and fixes them. Afterwards the solution is recomputed to optimize the remaining capacities. Before fixing a bike lane, it is ensured that the remaining car network remains strongly connected. Each iteration results in a feasible graph division into car and bike network, assuming all lanes aside from the fixed bike lanes are car lanes. 

\section{Experiments on real data}\label{sec:real_data}

\begin{figure}[htb]
    \centering
    \includegraphics[width=\textwidth]{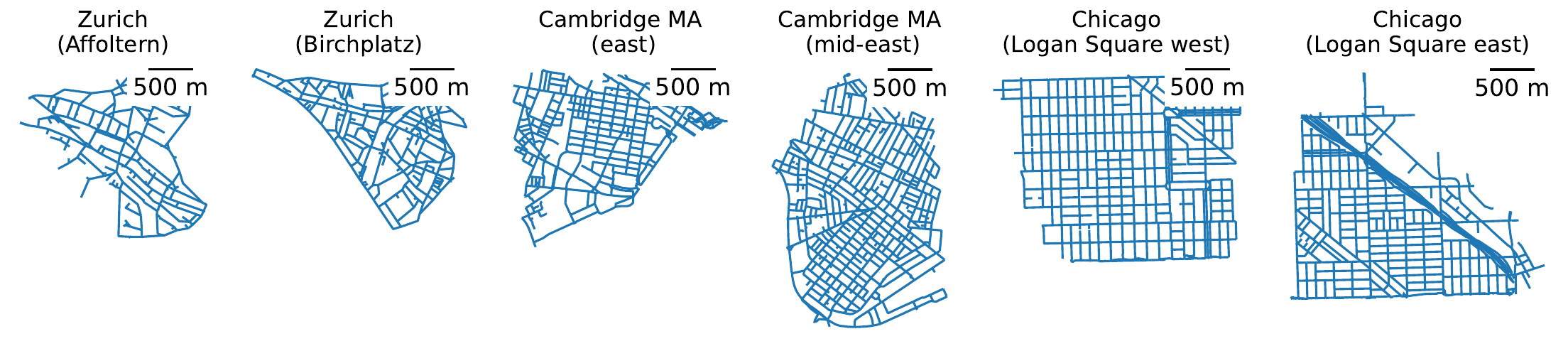}
    \caption{Network layout of the real instances used to test our algorithm}
    \label{fig:instance_graphs}
\end{figure}

We test the presented evaluation and optimization framework on real data from Zurich (Switzerland), Cambridge MA (US) and Chicago (US). From each city, two districts were selected, and their street network was extracted from Open Street Map (OSM) and pre-processed with the SNMan Python library\footnote{https://github.com/lukasballo/snman} by \citet{ballo2023} (see \autoref{sec:preprocessing}). An overview of the properties of the six instances and their graph layout is shown in \autoref{tab:instances} and \autoref{fig:instance_graphs}. Zurich, Affoltern, is a comparably small network with 529 edges, whereas the district Birchplatz already has 830 edges. For Cambridge MA, we consider one mid-sized instance and one large part with 1507 lanes. The instances in Chicago, Logan Square, differ from the other instances due to their strong grid-like layout (see \autoref{fig:instance_graphs}). The OD-matrices $\Omega$ are derived from public bike sharing data (Chicago \& Cambridge MA) or census data (Zurich); see \autoref{sec:preprocessing} for details. The full Pareto fontier is computed by executing \autoref{alg:rounding} once with $k=50$, i.e. re-solving the LP every 50 edge allocations. 

\begin{table}[htb]
\centering
\resizebox{0.8\textwidth}{!}{
\begin{tabular}{lcccccc}
\toprule
 & \multicolumn{2}{c}{Zurich}   & \multicolumn{2}{c}{Cambridge MA} & \multicolumn{2}{c}{Chicago}   \\
 & Affoltern & Birchplatz & east  & mid-east & L.S. west & L.S. east \\
\toprule
Nodes & 213 & 301 & 283 & 504 & 345 & 371 \\
Streets & 290 & 431 & 506 & 775 & 601 & 603 \\
Lanes & 535 & 799 & 883 & 1253 & 934 & 1126 \\
OD-paths & 219 & 383 & 456 & 677 & 175 & 281 \\
OD-paths extended & 430 & 684 & 738 & 1180 & 520 & 648 \\
Runtime optim. [min] & 4 & 36 & 20 & 146 & 26 & 38 \\
Runtime Pareto [h] & 0.28 & 4.3 & 2.74 & 7.91 & 5.35 & 7.07 \\
\bottomrule
\end{tabular}
}
\caption{Overview of real-world instances. The runtime for solving the LP one time is given as ``Runtime optimization'', in contrast to the runtime for computing the whole Pareto frontier (``Runtime Pareto'').}
\label{tab:instances}
\end{table}

In the following, we report the travel times over OD-paths in $\Omega$ if not denoted otherwise. All tests were executed on a standard machine with 16 GB RAM, using a Gurobi solver. The source code is available at \url{https://github.com/mie-lab/bike_lane_optimization}.

\newpage
\subsection{Heuristic baselines}\label{sec:heuristics}
\vspace{1em} 

To demonstrate our evaluation framework based on Pareto optimality, we compare our LP algorithm to three heuristic methods inspired by previous work. \citet{steinacker_demand-driven_2022} proposed to generate a sequence of bike networks, starting from a network with bike lanes at every edge, and iteratively removing bike lanes based on their betweenness centrality, an index measuring how frequently an edge is part of a shortest path~\citep{brandes2008variants}. In their work, the betweenness centrality is computed with respect to the shortest paths of an OD-matrix derived from the pickups and drop-offs in a bike sharing system. 
Our work builds up on this demand-driven approach; however, \citet{steinacker_demand-driven_2022} ignore the effect of new bike lanes on other traffic, which limits the comparability between their resulting networks and those produced by our algorithm. Despite this, their underlying method—betweenness centrality—can be adapted to develop an iterative algorithm that aims to allocate bike lanes in areas with minimal impact on the car network. We implemented three baseline approaches based on the betweenness centrality measure that iteratively allocate bike or car lanes. For the first one, we follow \citep{steinacker_demand-driven_2022} closely and start from a network where all streets are bike lanes. To express the negative impact of proper bike lanes on the car network, we exploit the concept of \textit{bike priority lanes}, where cars are required to give priority to bikes and are slowed down accordingly. Specifically, it is assumed that cars can drive 10 km/h on bike lanes. The initial network is thus a full bike lane network with a car speed of 10km/h throughout the city. As in \citep{steinacker_demand-driven_2022}, edges are iteratively removed from the network (and designated as proper car lanes), starting from the edge with lowest betweenness centrality with respect to the shortest bike travel times. The same OD-matrix and (perceived) travel times as for our optimization approach are used to ensure comparability. Since this method starts from a full bike network, we call this approach \textit{betweenness-top-down}.

In contrast, for the second and third baseline, we start from a network with car lanes only and add bike lanes iteratively, termed \textit{betweenness-bottom-up} in the following. There are two ways to re-assign lanes to cycling: 1) a car-prioritizing approach, where the first edges to be converted to bike lanes are the ones with \textit{lowest} betweenness centrality, computed with respect to the \textit{car} travel time (\textit{betweenness-bottom-up (car)}), and 2) a bike-prioritizing approach, where the edges with \textit{highest} betweenness centrality with respect to the \textit{bike} travel time are converted first (\textit{betweenness-bottom-up (bike)}). In both cases, we iteratively select the edge with the lowest/highest betweenness centrality computed on the OD-matrix, and convert this edge into a bidirectional bike lane, setting the car travel time along this edge to $\infty$ as in our optimization approach. In summary, we compare our approach to three strong heuristic methods based on the betweenness centrality, where the first one (top-down) is designed to be as similar as possible to \citep{steinacker_demand-driven_2022} while the second and third one (bottom-up) are constructed for better comparability to our method. The method \textit{betweenness-bottom-up (car)} resembles the approach taken in \cite{ballo2023}. Flow charts outlining the procedure for each baseline are provided in Appendix~\ref{sec:baselines_app}.

\vspace{1em} \subsection{Model comparison by Pareto optimality}\label{sec:results_real}
\vspace{1em}

\begin{figure}[htb]
    \centering
    \includegraphics[width=0.9\textwidth]{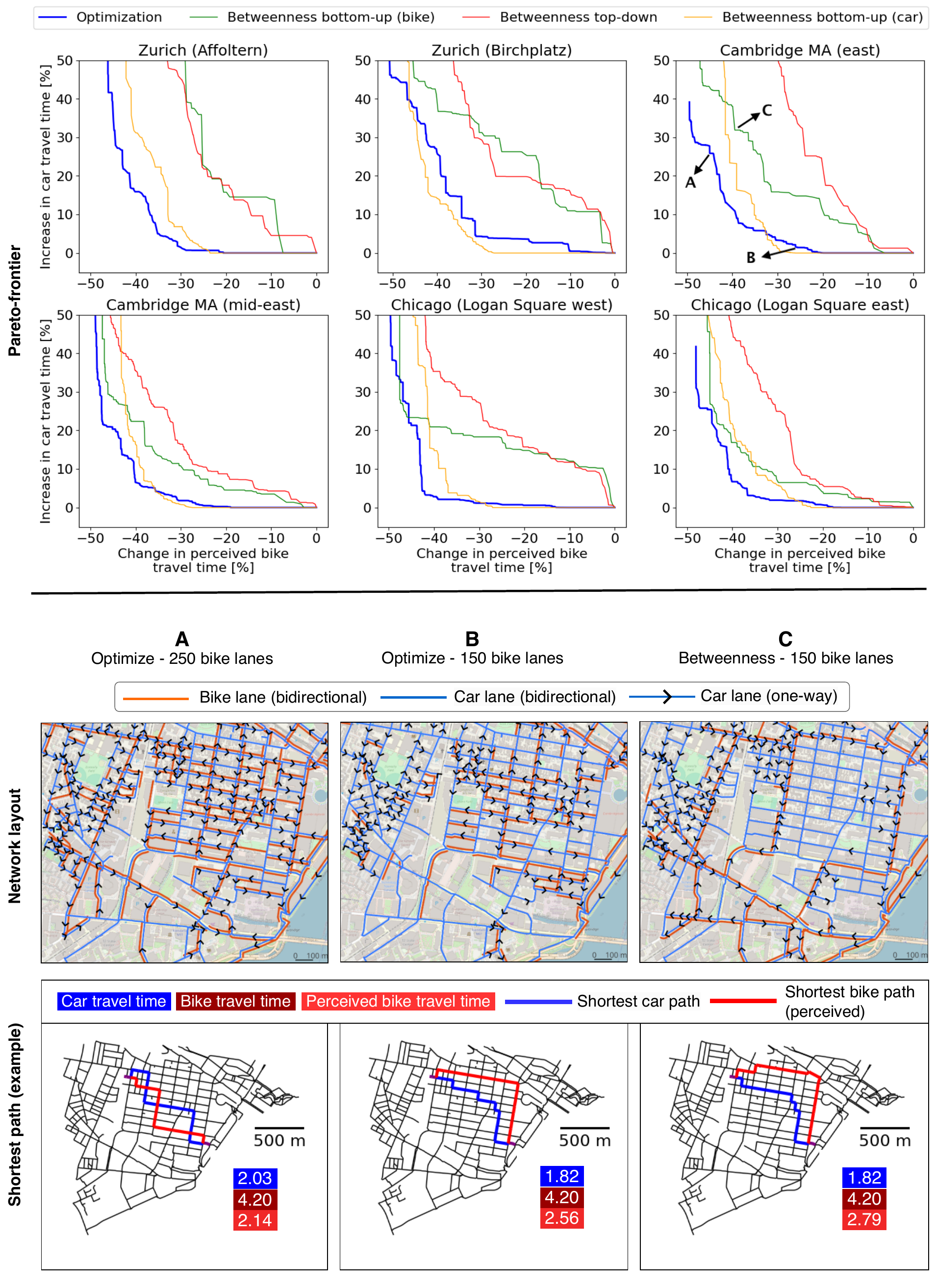}
    \caption{Pareto optimality of bike networks. Top: Algorithms are compared by their Pareto frontier. In five out of six instances, our linear programming approach outperforms methods based on the betweenness centrality. Bottom: Each point on the Pareto frontiers (top) corresponds to one plausible street network. Three examples in Cambridge MA are shown, where the bike networks differ dependent on the planning method and the number of allocated bike lanes. This is also reflected in the distance of shortest paths, where the existence of dedicated cycling infrastructure is rewarded in the perceived bike travel time.}
    \label{fig:real_world}
\end{figure}

\autoref{fig:real_world} presents the Pareto frontier for each method and instance, visualizing the achieved trade-off between bike and car travel time. Each point on the Pareto frontier is one bike network. 
For example, for Zurich Affoltern our algorithm yields a network where the perceived bike travel time is decreased by 40\% while increasing the car travel time only by 17\%. It is worth noting that the bike travel time can be reduced by at most 50\% due to the setting of $t^{\beta}(e) = 2 t^b(e)$. In five instances, the Pareto frontier of our optimization approach dominates over the heuristic solutions. From the baselines, the betweenness-top-down approach performs worst. This approach strongly depends on the allowed car speed on bike-priority lanes, which is only set to 10km/h in our experiment. The bike-focused bottom-up method usually yields better networks when many bike lanes are allocated, whereas the car-focused bottom-up method yields better solutions when few bike lanes are allocated (see intersections of green and yellow lines for Cambridge MA and Chicago). The framework thus also provides recommendations about the use cases for the respective methods.

\autoref{fig:real_world} illustrates three distinct bike networks proposed for Cambridge MA (east), varying due to the planning method and the number of allocated bike lanes. The betweenness-centrality method (C) places a greater number of bike lanes on main roads compared to our optimization approach (B). When allocating more bike lanes (250 instead of 150) with the optimization approach, the car network transitions into a complex one-way system. This change leads to increased car travel times, as depicted with one example in \autoref{fig:real_world}, with travel times rising from 1.82 to 2.03 for a specific route. Although the overall bike travel time remains consistent across all algorithms—since cyclists have access to all roads—the perceived bike travel time varies significantly. Specifically, the perceived travel time is considerably longer when using the betweenness algorithm compared to the optimization method, highlighting the impact of planning approaches on travel efficiency.

\FloatBarrier

\vspace{1em} \subsection{Analyzing the bike-car trade-off between cities}

\begin{figure}[htb]
    \centering
    \includegraphics[width=0.6\textwidth]{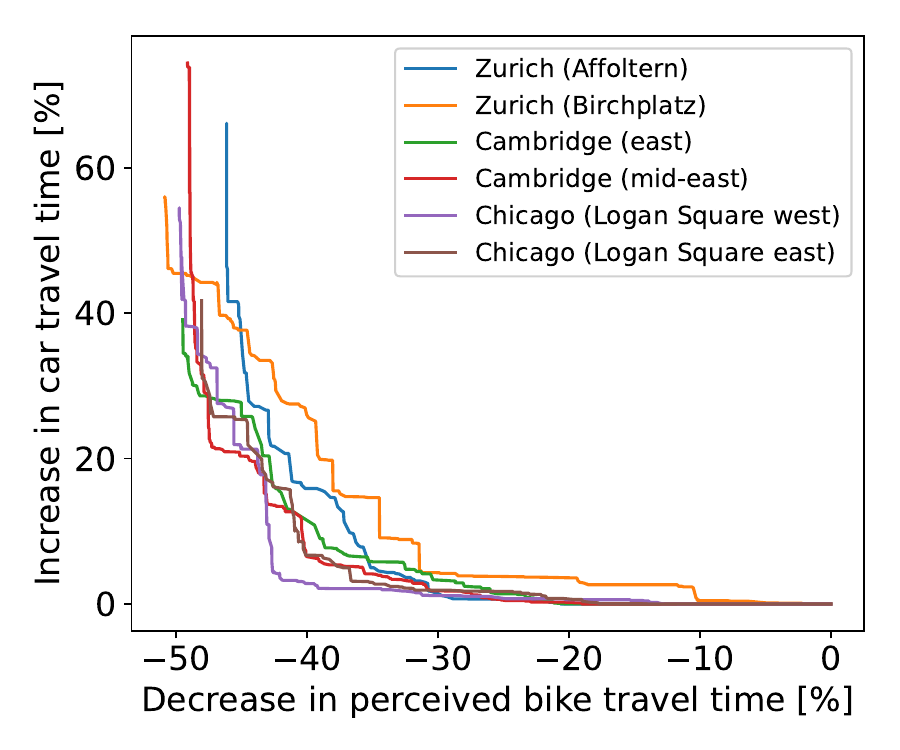}
    \caption{Comparison of the trade-off between bike and car travel times in different cities.}
    \label{fig:real_world_comp}
\end{figure}

The evaluation framework also assesses the adaptability of urban areas to bike lane integration, highlighting the variability in the ease of adding bike lanes without significantly disrupting car traffic. Figure \ref{fig:real_world_comp} showcases a comparison of the six instances through their Pareto frontiers with respect to the optimization algorithm. For instance, in Chicago Logan Square (west), it is possible to reduce perceived bike travel time by 48\% with only a 12\% increase in theoretical car travel time. Conversely, in Zurich-Affoltern and Cambridge MA-East, integrating bike lanes proves considerably more challenging. The Pareto frontier for several cities reveals points at which car travel time rises sharply, providing valuable insights for urban planners on the optimal number of bike lanes to introduce. Across all tested districts, improving cycling conditions by more than 35\% is achievable without exceeding a 20\% increase in car travel time.

An interesting question for future work is to understand the \textit{properties} of a street network that impact the shape of its Pareto frontier. As a simple attempt towards this goal, we computed several indices for network characteristics, 
selected from a list of metrics provided by \citet{boeing2022street}. The results are provided in \autoref{tab:data_analysis}, where the solution quality is measured in terms of three intersections of the Pareto frontier, namely the car travel time increase at a bike time decrease of 35, 40 and 45\%. Due to the small sample size of six instances, the statistical significance is limited, but, for example, we observe a relation of the solution quality to the betweenness centrality, favoring the use of this heuristic. There is also a small negative correlation between the average node degree and the goodness of the Pareto frontier. In general, the differences indicate that US cities are easier to optimize than Zurich as a European example. 

\begin{table}[htb]
    \centering
\resizebox{\textwidth}{!}{
\begin{tabular}{ll|rrrrrr}
\toprule
 &  & Chicago & Cambridge MA & Chicago & Cambridge MA & Zurich & Zurich \\
 & & (L.S. west) & (mid-east) & (L.S. east) & (east) & (Affoltern) & (Birchplatz)\\
\toprule
\midrule
\multirow[t]{2}{*}{Node degree} & Mean & 3.316 & 2.968 & 3.008 & 3.004 & 2.364 & 2.691 \\
 & StD & 0.895 & 0.921 & 0.945 & 0.950 & 1.093 & 1.025 \\
\midrule
\multirow[t]{2}{*}{Gradient} & Mean & 1.849 & 0.000 & 2.111 & 0.000 & 2.724 & 2.261 \\
 & StD & 2.215 & 0.000 & 2.498 & 0.000 & 3.055 & 2.248 \\
\midrule
\multirow[t]{2}{*}{Street width [m]} & Mean & 4.777 & 5.522 & 6.307 & 6.684 & 5.307 & 5.696 \\
 & StD & 2.475 & 2.894 & 4.413 & 4.692 & 2.671 & 2.959 \\
\midrule
\multirow[t]{2}{*}{Street length [km]} & Mean & 0.129 & 0.095 & 0.129 & 0.099 & 0.116 & 0.088 \\
 & StD & 0.066 & 0.062 & 0.085 & 0.061 & 0.101 & 0.068 \\
\midrule
\multirow[t]{2}{*}{Closeness} & Mean & 0.088 & 0.078 & 0.088 & 0.089 & 0.098 & 0.094 \\
 & StD & 0.014 & 0.010 & 0.012 & 0.013 & 0.014 & 0.014 \\
\midrule
\multirow[t]{2}{*}{Betweenness} & Mean & 0.020 & 0.017 & 0.021 & 0.027 & 0.041 & 0.027 \\
 & StD & 0.021 & 0.019 & 0.022 & 0.027 & 0.042 & 0.031 \\
\midrule
Clustering & Mean & 0.023 & 0.021 & 0.029 & 0.035 & 0.032 & 0.021 \\
coefficient & StD & 0.077 & 0.079 & 0.084 & 0.110 & 0.111 & 0.076 \\
\midrule \midrule
\textbf{Car time increase (in\%)}& -35\% & 2.124 & 4.150 & 2.985 & 5.803 & 4.988 & 14.647 \\
\textbf{corresponding to} & -40\% & 2.820 & 6.515 & 6.715 & 10.857 & 15.890 & 26.049 \\
\textbf{bike time decrease of...}& -45\% & 21.310 & 20.404 & 25.406 & 25.831 & 34.399 & 37.673 \\
\bottomrule
\end{tabular}
}
    \caption{Analysis of street network properties with respect to the Pareto frontier. Indicators were selected based on \citet{boeing2022street}.}
    \label{tab:data_analysis}
\end{table}
\FloatBarrier

\subsection{Case study: Rebuilding a whole city}\label{sec:case_study}

\vspace{1em}

One limitation of the proposed algorithm is its runtime complexity, hindering its application on metropolitan street networks of ten thousands of edges. However, the city can be divided into districts and processed region-wise, allowing parallelized optimization. As a case study, we aim to allocate as many bike lanes as possible in the whole city of Zurich, a street network of 7591 nodes and 9949 edges in total, resulting in a lane graph of 17870 edges. To augment the Pareto optimality analysis provided for the six real-world instances above, the aim of this case study is, on the one hand, to show how divide-and-conquer strategies make the algorithm applicable for large-scale networks and, on the other hand, to explore the extreme scenario of transforming Zurich into a predominantly (e-)bike city, as discussed in recent work~\citep{ballo_e-bike_2023}. 

We process the main roads separately from non-main roads, where the city is partitioned into five parts and into 57 regions respectively (for details see \autoref{app:partitioning}). For each region, \autoref{alg:rounding} is applied until convergence; i.e., the LP is solved and bike lanes are allocated in several iterations until any further bike lane allocation would disconnect the network. The algorithm execution is parallelized with the Python \texttt{multiprocessing} library and is executed on five cores and took 16h 35min in total. While most regions are relatively small (averaging 402 edges and 191 OD pairs), some partitions are still comparable in size to the largest instances tested in \autoref{sec:results_real}. For example, one main-road network contains 1470 edges and 1492 OD pairs. Thus, considering that 62 instances must be solved with a total size of 17870 edges, the total runtime is surprisingly low.

\begin{figure}[htb]
    \centering
    \begin{subfigure}[b]{0.49\textwidth}  
    \includegraphics[width=\textwidth]{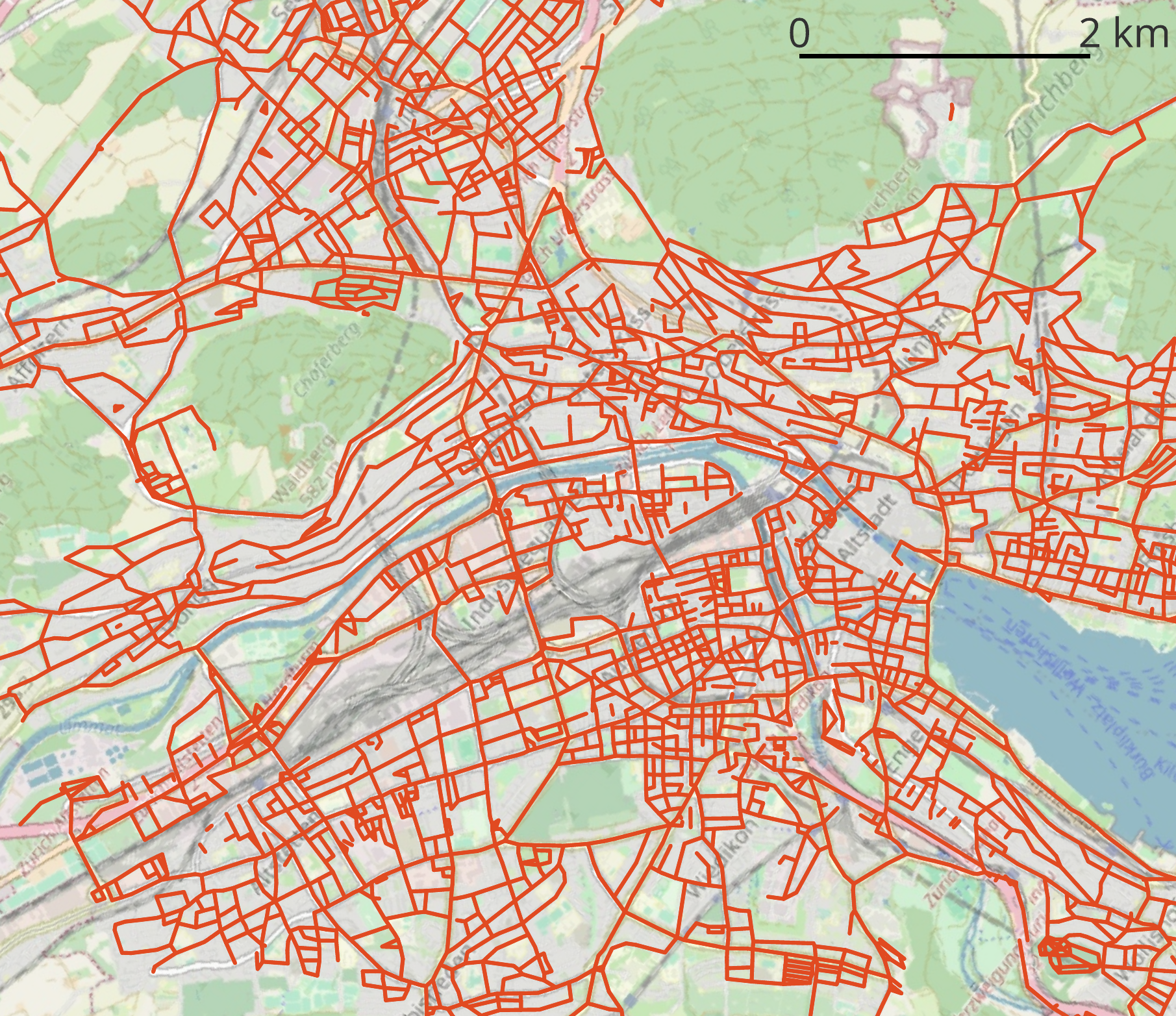}
    \caption{Bike lane network (bidirectional)}
    \label{fig:zurich_bikes}
    \end{subfigure}
    \begin{subfigure}[b]{0.49\textwidth}  
    \includegraphics[width=\textwidth]{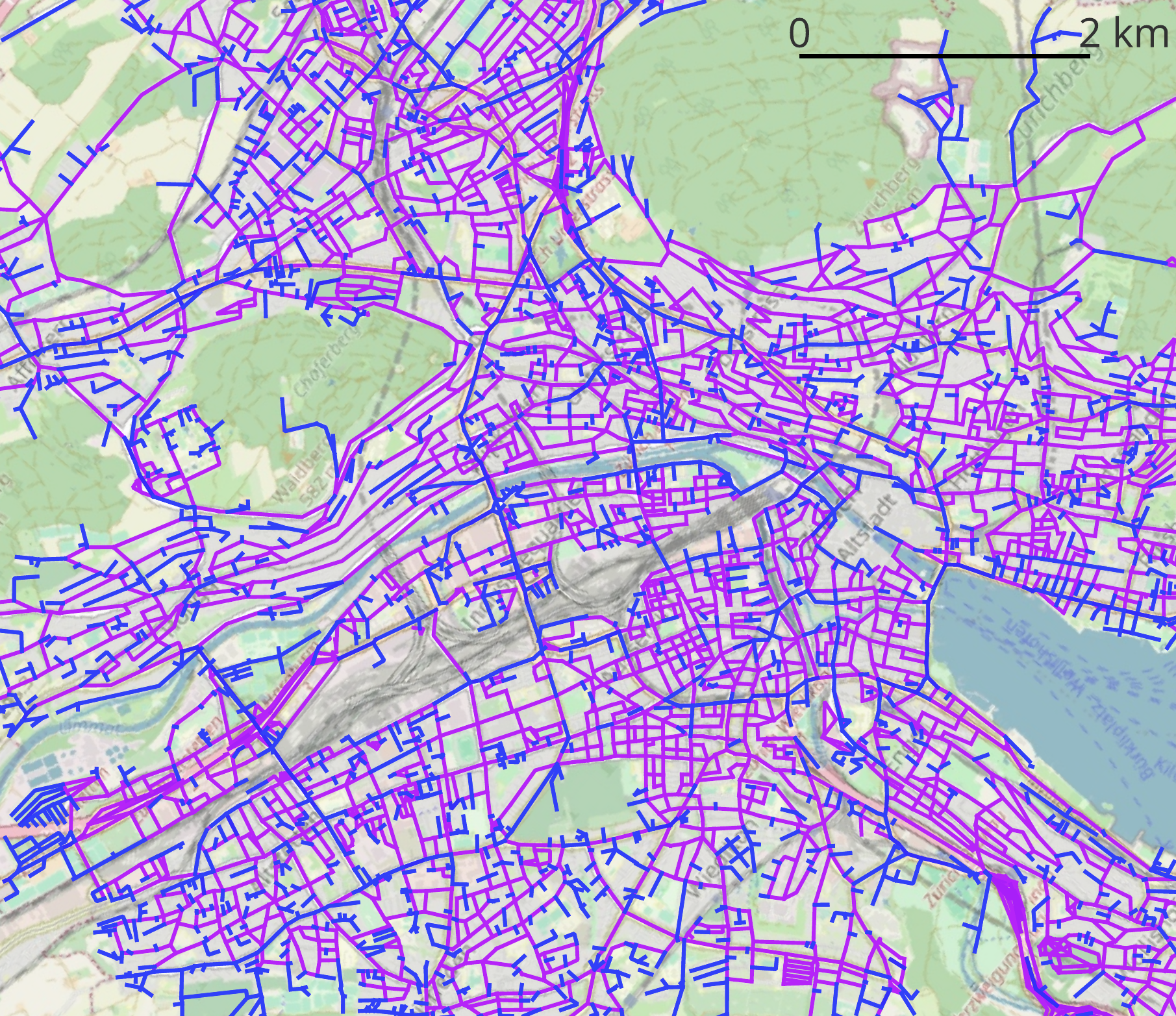}
    \caption{Car lane network (pink: one-way, blue: bidirectional)}
    \label{fig:zurich_cars}
    \end{subfigure}
    \caption{Case study of rebuilding the whole city of Zurich. As many bike lanes as possible were allocated with our optimization approach. The dense bike network results in many one-way car lanes, nevertheless forming a functional and connected car network.}
    \label{fig:zurich_network}
\end{figure}

The transformed street network is shown in \autoref{fig:zurich_network}. Allocating as many bike lanes as possible in each district results in a dense and extensive bike network (\autoref{fig:zurich_bikes}), at the cost of a car network composed of many one-way streets (\autoref{fig:zurich_cars}). After applying the algorithm, 60\% of the streets are one-way streets for cars, compared to 15\% before. To understand the effect of such a radical re-distribution of space, we analyze the network change in terms of the added space and the change in (perceived) travel time. \autoref{tab:metrics_case_study} shows that our algorithm is able to allocate 559km bike lanes while ensuring that the car network remains connected. Highways and trams are fixed, but around one third of the motorized roads is repurposed for bike lanes, with a similar ratio for main and local roads. Crucially, this decreases the perceived bike travel time by 43\%, while the car travel time increases by 39\%. The potential for bike lane allocation varies by region, with an average of 33\% of total road kilometers allocated to bike lanes (standard deviation of 9\%) after applying the algorithm. In one region, 54\% of the road space is allocated to bike lanes. Together with the results provided in \autoref{sec:results_real}, this demonstrates that the success of bike lane allocation using our algorithm—and likely with any algorithm—depends heavily on the specific network layout and its properties. Future work could further quantify cities' potential to transform into highly bike-friendly urban areas by applying our algorithm across different urban networks.

\begin{table}[htb]
    \centering
\resizebox{0.8\textwidth}{!}{
\begin{tabular}{l|ccccc|cc}
\toprule
& \multicolumn{5}{c}{\textbf{Available space (in km)}} & \multicolumn{2}{c}{\textbf{Avg. travel time [min] }} \\ 
 & Cycling & Motorized & Motorized  & Highways & Tram & Bike & Car \\
 &  &   main roads & local roads  &  &  & \textit{perceived}  & \\
\midrule
\textbf{Before} & 17.03 & 513.79 & 1171.16 & 84.61 & 107.55 & 23.47 & 5.22 \\
\textbf{After} & 559.15 & 349.72 & 782.85 & 84.61 & 107.55 & 13.48 & 7.29 \\
\bottomrule
\end{tabular}
}
\caption{Metrics for radical space distribution of a whole city}
\label{tab:metrics_case_study}
\end{table}

\FloatBarrier

\vspace{2em} \section{Experiments on synthetic data}\label{sec:synthetic}

In the following, we will present experiments demonstrating the computational efficiency and solution quality of the proposed algorithm and relaxations. For the sake of a systematic evaluation, we generate synthetic street networks. 
Given a desired number of nodes $n$, we sample $n$ pairs of spatial coordinates uniformly. To resemble real street networks in a simple fashion, we sample 2, 3, or 4 neighbors (count decided randomly) for each node $u$, where the sample probability of a node $v$ is proportional to its inverse quadratic distance to $u$, formally $P(v) = \frac{d(u, v)^{-2}}{\sum_v d(u, v)^{-2}}$, except for $P(u)$ which is set to zero to avoid loops. Whenever a node $v$ is sampled as a neighbor of $u$, a street with two lanes is added between them, corresponding to an edge of capacity 2. The neighbors are sampled with replacement to allow multi-lane streets between nodes. If the generated graph is not strongly connected, the sampling is repeated. Each street is assigned a distance and gradient based on the spatial coordinates. 

\vspace{1em} \subsection{Empirical runtime}
\vspace{1em}

Due to the constraints in \autoref{eq:flowbike} and \autoref{eq:flowcar}, there are $3 m \cdot |\Omega|$ flow variables to be optimized, in addition to $2m$ capacity values, where $m=|E'|$ is the number of edges. To analyze the empirical runtime of our approach, we generate synthetic networks according to the scheme described above and systematically vary the number of variables. Let $\nu$ be the total number of variables, with $\nu = 3 \cdot m \cdot |\Omega| + 2 \cdot m$. We set the number of $\Omega_{demand}$ in relation to its maximum size $n^2$, such that $|\Omega_{demand}| = \psi \cdot n^2$ with $\psi \in [0,1]$. However, $n$ paths are further added with $\Omega_{aux}$ to ensure strong connectivity (see \autoref{sec:methods_od_matrix}). Thus, dependent on $\psi$, the total number of variables amounts to $\nu = 3 \cdot m \cdot |\Omega| + 2 \cdot m = 3 \cdot m \cdot (\psi n^2 + n) + 2 \cdot m$. For our experiment, we sample $\nu\in \{0.4M, 0.7M, 1M, 1.3M, 1.6M, 1.9M, 2.2M , 2.5M, 2.8M\}$ (with $M$ = million) and for each of them generate random networks with $n\in \{100, 150, 200, 250, 300, 350, 400\}$. Solving the above equation for $\psi$, we set $\psi = \frac{\nu - 2 * m}{n^2 * m * 3} - \frac{1}{n}$ dependent on the network size and the desired number of variables, and run the optimization with a randomly generated OD-matrix $\Omega$ of the respective size. 

\begin{figure}[htb]
    \centering
    \begin{subfigure}[b]{0.48\textwidth}
    \includegraphics[width=\textwidth]{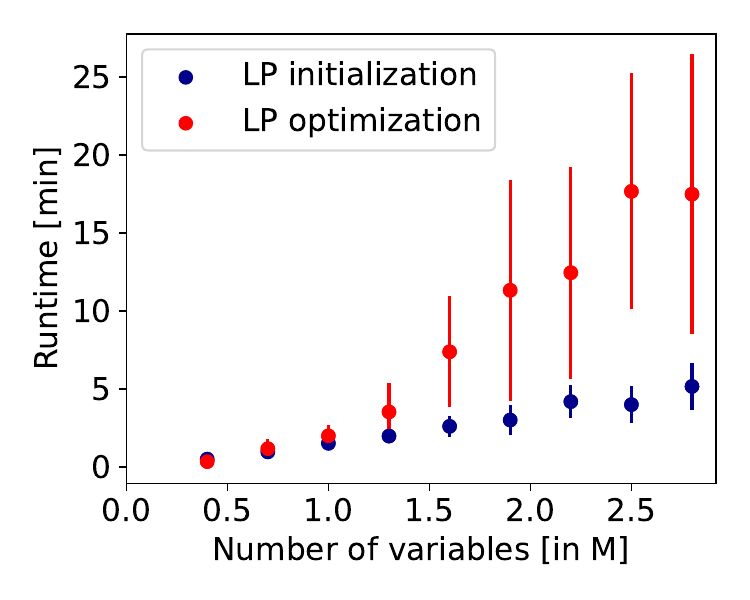}
    % \caption{By number of OD-pairs}
    % \label{fig:od_runtime}
    \caption{Mean and StD of LP initialization and optimization}
    \label{fig:runtime_init_optim}
    \end{subfigure}
    \hfill
    \begin{subfigure}[b]{0.48\textwidth}
    \includegraphics[width=\textwidth]{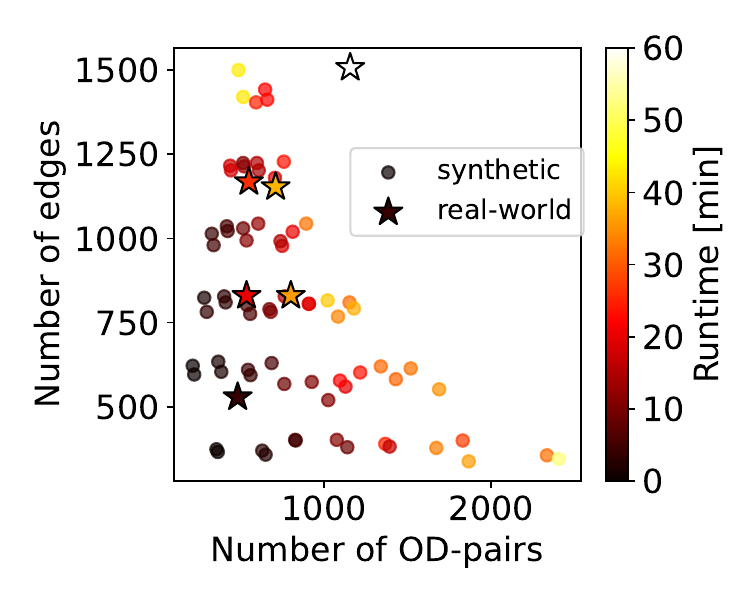}
    % \caption{By number of edges}
    % \label{fig:edges_runtime}    
    \caption{Runtime by the number of edges and s-t paths}
    \label{fig:runtime_edges_od}
    \end{subfigure}
    \caption{Empirical runtime analysis}
    \label{fig:runtime}
\end{figure}

In \autoref{fig:runtime_init_optim}, the corresponding runtimes are shown by the number of variables. The time for initializing the LP scales approximately linearly with the number of variables, as expected. 
Solving the linear program requires superlinear runtime with respect to the number of variables, in alignment with theoretical considerations of linear programming~\citep{cohen2021solving}. 
Nevertheless, a graph with 300 nodes, 1200 edges, and 450 node-pairs in $\Omega$ can still be optimized within ten minutes. % 30 and 82 seconds
\autoref{fig:runtime_edges_od} gives an intuition on the evolvement of the total runtime (LP initialization plus optimization) dependent on edge count and size of $\Omega$. Networks with a large number of edges can be processed efficiently when $\Omega$ is sufficiently small. The runtime of real-world instances aligns with the runtime for synthetic instances of similar size. 
Considering a typical size of 500 -- 1000 edges for city districts (see \autoref{tab:instances}), the experiments confirm the applicability of our algorithm to real-world instances. The corresponding Integer Problem, in contrast, already requires 1-2h to solve with instances of 500k variables (see comparison in Appendix~\ref{sec:integer_runtime}).

\FloatBarrier

\vspace{1em} \subsection{Effect of linear relaxation and rounding}\label{sec:synthetic_lp_vs_ip}
\vspace{1em}

The proposed algorithm is designed to optimize the trade-off between the goodness of the car and bike networks. This is expressed in two ways: the parameter $\gamma$ that defines the importance of the car travel time, and the number of edges actually allocated for bikes by the iterative post-processing algorithm (\Cref{alg:rounding}). With an \textit{integer} problem formulation, only $\gamma$ can be varied since post-processing is unnecessary. In the following, we aim to compare the solutions of our LP-based algorithm with the optimal output of an IP.

\begin{figure}
    \centering
    \begin{subfigure}[b]{0.49\textwidth}
    \includegraphics[width=\textwidth]{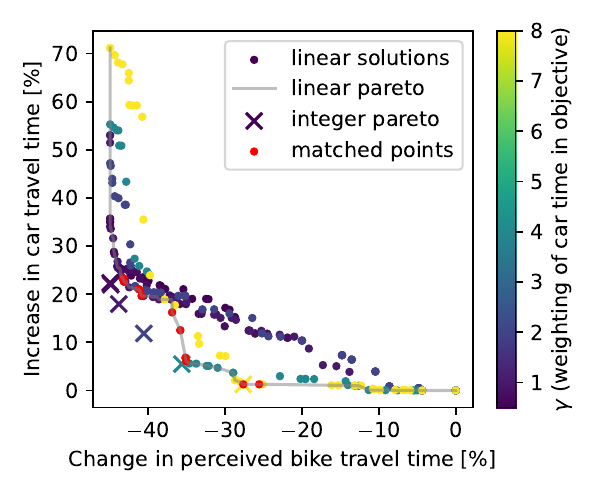}
    \caption{Example of the Pareto frontiers from linear and integer formulations}
    \label{fig:Pareto_int_lin_example}
    \end{subfigure}
    \hfill
    \begin{subfigure}[b]{0.49\textwidth}
    \includegraphics[width=\textwidth]{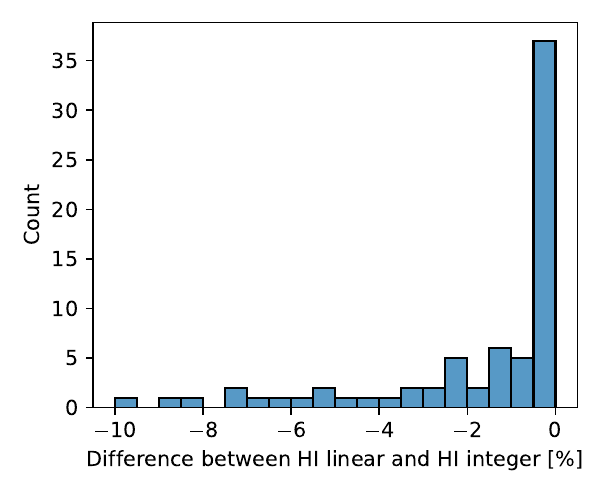}
    \caption{Distribution of percentage change from integer HI to the HI of the corresponding linear Pareto frontier.}
    \label{fig:hi_dist}
    \end{subfigure}
    \caption{There are two potential variations leading to Pareto frontiers: 1) the weighting in the objective functions and 2) the number of edges added in the rounding algorithm. \subref{fig:Pareto_int_lin_example} provides an example. The closest linear solutions (red) are compared to the Pareto frontier of the integer solver in terms of the hypervolume indicator (HI). The distribution of the difference in HI is shown in \subref{fig:hi_dist}. In most cases, the integer solution is only marginally better than the linear solutions, increasing the IP HI by less than 5\% for all tested instances.}
    \label{fig:lin_vs_int}
\end{figure}

\autoref{fig:Pareto_int_lin_example} provides an example of the outputs of the LP and the IP for one synthetic instance. Each point is one solution, i.e. one bike network. Solving the IP results in a set of Pareto-optimal solutions due to varying $\gamma$ (crosses in \autoref{fig:Pareto_int_lin_example}). The higher $\gamma$, the lower the average car travel time and the higher the bike travel time. In comparison, running the LP formulation with a specific $\gamma$ does not yield a single solution, but a whole Pareto frontier since bike lanes are allocated iteratively (points in \autoref{fig:Pareto_int_lin_example}). 
By selecting all linear solutions that are not dominated, a single Pareto frontier for the LP is constructed (gray line in \autoref{fig:Pareto_int_lin_example}). 

The Pareto frontiers for IP and LP are compared quantitatively via the hypervolume indicator (HI). The HI measures the area that is dominated by an n-dimensional Pareto curve with respect to a reference point. Let $(x_{ref}, y_{ref})$ be a two-dimensional reference point, here a tuple of bike and car travel time. The set of dominated solutions $A$ with respect to a Pareto frontier $P$ is defined as follows: $$A = \{(x, y)\ |\  x \leq x_{ref}, y \leq y_{ref},  \exists (\hat{x},\hat{y})\in P \text{ with } (\hat{x},\hat{y}) \text{ dominating } (x, y) \}$$
The HI is defined as the area comprised by $A$, spanning the area between the Pareto frontier (with rectangular interpolation) and the reference point. The higher the HI, the better the algorithm. In the following, we set the reference point to $(x_{ref}, y_{ref}) = (0, 50)$, i.e. to 0\% bike travel time decrease and 50\% car travel time increase. Although the algorithm also produces solutions with more than 50\% car travel time increase, the shape of the Pareto frontier below 50\% is arguably more relevant for urban planning. 
Problematically, a Pareto frontier with more points yields a higher HI than a frontier with few points. By design of the experiment, the Pareto frontier for the LP is more dense (see \autoref{fig:Pareto_int_lin_example}). For a fair comparison, we thus select only a subset of Pareto-optimal linear solutions, namely the ones that are closest to the integer solutions (red points in \autoref{fig:Pareto_int_lin_example}).

To compare the linear and the integer solutions systematically, the HI is compared on 60 synthetic instances with $n\in[20, 40, 60, 80]$ nodes (due to the large runtime of the integer problem, only small graphs are tested), with 10\% of all possible OD-pairs for $n\in[20, 40]$ and 1\% for $n\in[60, 80]$. The post-processing algorithm (\Cref{alg:rounding}) for rounding fractional solutions of the LP is run with $k=10$, i.e., the LP is re-optimized every 10 edges allocated as bike lanes. The HI is computed for the integer Pareto frontier and the selected linear points, and their percentage difference is computed as $(HI_{linear} - HI_{integer}) / HI_{integer} * 100$. The distribution of the difference is given in \autoref{fig:hi_dist}. In 78\% of the tested instances, the linear HI is less than 3\% smaller than the integer HI. In 23\% cases, the integer solutions are fully recovered. Thus, the linear relaxation offers a good trade-off between efficiency and optimality, yielding networks that are almost on-par with the integer solution in terms of bike and car travel times. In Appendix~\ref{app:int_vs_lin} we provide additional analysis to clarify the conditions under which the optimal solution can be recovered. Our findings indicate that the size of $\Omega$—the number of OD pairs—plays a critical role. As the number of OD pairs increases, the task of allocating bike and car lanes becomes more complex, leading to greater divergence between the LP and IP solutions.

\FloatBarrier

\vspace{1em} \subsection{Baseline comparison}\label{sec:synth_baselines}
\vspace{1em}

\begin{figure}
    \centering
    \begin{subfigure}[b]{0.49\textwidth}
    \includegraphics[width=\textwidth]{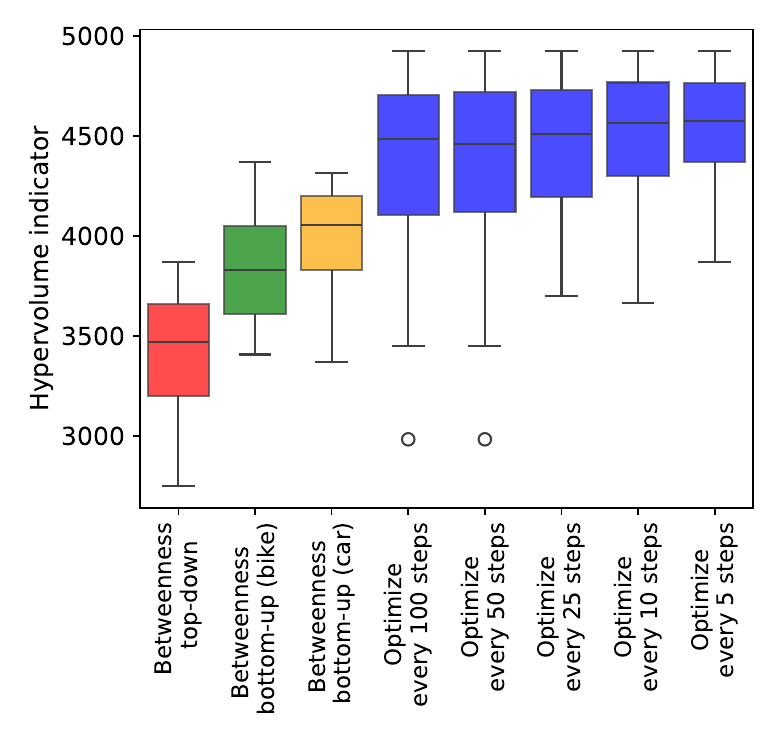}   
    \caption{Betweenness vs Optimization}
    \label{fig:compare_rounding_algorithm}
    \end{subfigure}
    \begin{subfigure}[b]{0.49\textwidth}
    \includegraphics[width=\textwidth]{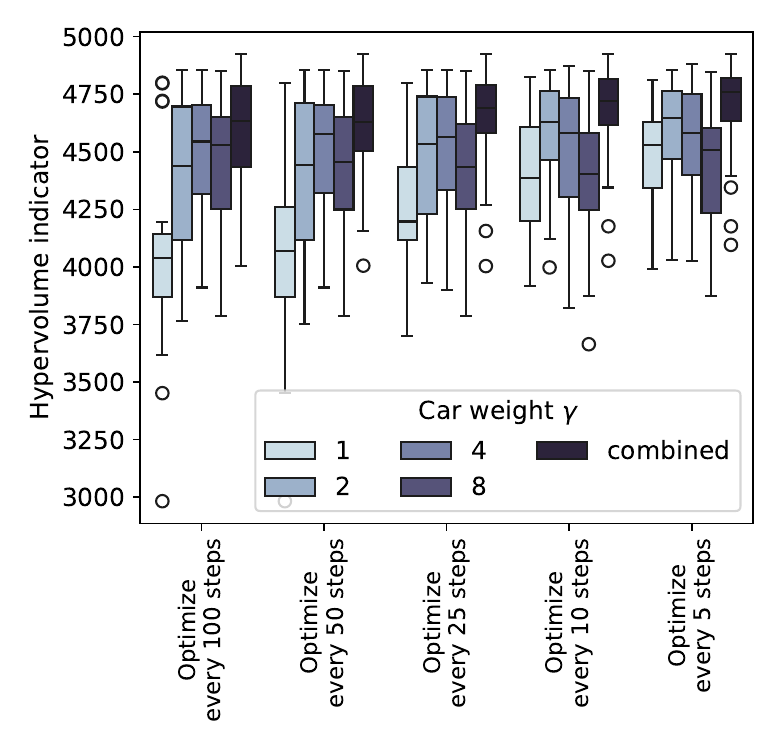}   
    \caption{Optimization by car weight}
    \label{fig:compare_rounding_carweight}
    \end{subfigure}
    \caption{Algorithm comparison by hypervolume indicator (HI). In (\subref{fig:compare_rounding_algorithm}) we compare the optimization approach (blue) to three baselines based on the betweenness centrality. The optimization approach generally performs better than the baselines, and re-optimizing more frequently improves performance. The best results are achieved when combining the results for different travel time weighting determined by the $\gamma$ parameter (\subref{fig:compare_rounding_carweight}).}
    \label{fig:alg_comp}
\end{figure}

Furthermore, we benchmark our algorithm against three heuristic algorithms introduced in \autoref{sec:heuristics}, using 20 synthetic instances ($n\in \{20, 30, 40, 50\}$, with five instances per $n$). For each instance, $\frac{1}{10}n^2$ OD pairs are sampled randomly to construct $\Omega$. \autoref{fig:compare_rounding_algorithm} compares the algorithms by the distribution of their hypervolume indicator. Each point in the plot represents one Pareto frontier. Our method's results, varying by the parameter $k$ (which influences the frequency of re-solving the linear program as described in Algorithm \ref{alg:rounding}), show a clear advantage over the baselines. As anticipated, a lower $k$ value in our approach yields better outcomes. Notably, the car-centric baseline algorithm performs best among the heuristic approaches, in accordance with the results on real data. 

\autoref{fig:compare_rounding_carweight} distinguishes our results by the car weight $\gamma$. Setting $\gamma$ between 2 and 4 yields good performance, whereas setting $\gamma=1$ or $\gamma=8$ shifts the Pareto frontier and fails to find the Pareto-optimal solutions that solving the IP would yield. The best solutions are obtained when combining the solutions from several car weights into one Pareto frontier, as shown in \autoref{fig:Pareto_int_lin_example} (gray line). Interestingly, combining the runs with different car weights for $k=100$ matches the performance of individual runs with $k=5$. However, it depends on the graph size whether it is more efficient to re-optimize more often (i.e., lower $k$) or to execute the whole process with several car weightings $\gamma$, since both requires to solve the LP multiple times. 

\vspace{1em} \subsection{Effect of demand-driven shortest path computation}
\vspace{1em}

The reduction from the all-pairs-shortest-path problem to an OD-shortest path problem via $\Omega$ can be viewed as a necessary relaxation with respect to the runtime, or as a feature of the algorithm to incorporate real-world travel demand. While we argue for the latter, it is interesting to analyze the effect of this relaxation on the all-pairs travel length. Are the omitted s-t-paths optimized implicitly, or do their travel times increase tremendously? To answer this question, we optimize the LP with OD matrices of varying size for random instances, and, for each of them, compute the objective value with respect to the all-pairs-shortest-path. This objective value is strictly larger when considering only a subset of the possible s-t-paths. Equality can only be achieved if the optimal capacity values for edges along omitted s-t-paths by chance correspond to the ones that minimize the all-pairs shortest path objective. 

\begin{figure}[htb]
    \centering
    \begin{subfigure}[b]{0.59\textwidth}
        \includegraphics[width=\textwidth]{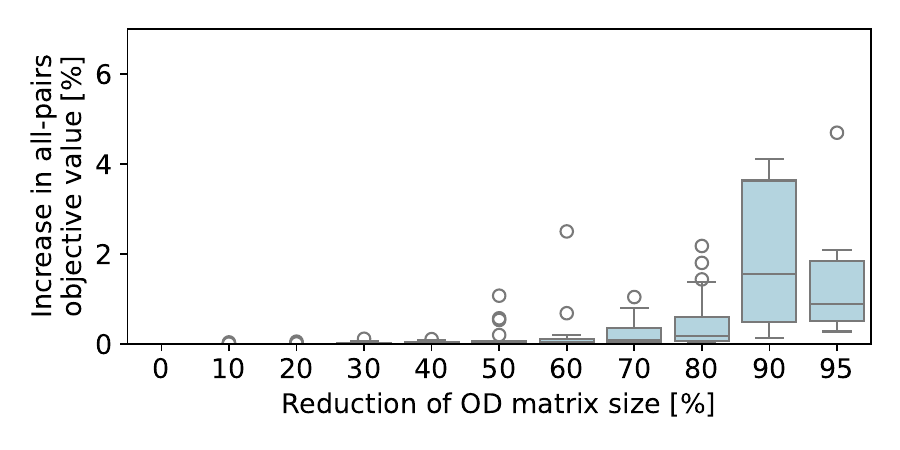}
        \caption{By number of s-t-paths}
        \label{fig:od_dependency_od}
    \end{subfigure}
    \hfill
    \begin{subfigure}[b]{0.39\textwidth}
        \includegraphics[width=\textwidth]{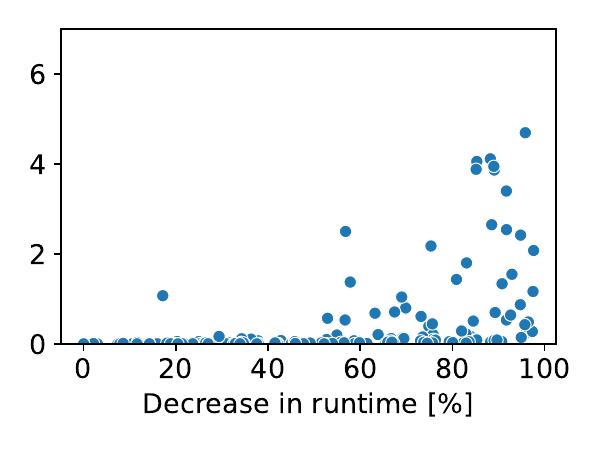}
        \caption{By runtime}
        \label{fig:od_dependency_runtime}
    \end{subfigure}
    \caption{Increase in the all-pairs-shortest-paths objective value, in relation to the decrease of OD-matrix size or runtime.}
    \label{fig:od_dependency}
\end{figure}

In \autoref{fig:od_dependency}, the increase in the objective value is contrasted to the decrease in OD-matrix size, and, resulting from that, the decrease in runtime. Each point corresponds to one random graph instance ($n\in \{20, 30, 40, 50\}$) together with an OD-matrix. Specifically, we tested 16 graphs (four of each $n$), and for each of them we set the size of the OD-matrix to 5\%, 10\%, 20\%, ..., 90\%, 100\% of the all-pairs matrix respectively\footnote{The exact size of each OD-matrix varies slightly due to the insertion of auxiliary OD-pairs to ensure connectivity (see \autoref{sec:methods_od_matrix}). 
}.

\autoref{fig:od_dependency} shows that the all-pairs travel times only increase marginally for reduced OD matrices. Even with up to 60\% reduction in the number of considered s-t-pairs or in the runtime, the objective value increases by less than 1\% for a majority of trials. Since the runtime scales linearly with the size of the OD-matrix, there is a similar relation between runtime and objective value (\autoref{fig:od_dependency_runtime}). The results indicate that the s-t paths in the reduced OD matrices still cover a significant portion of other potential paths in the graph, whose respective travel times are implicitly minimized. For real-world instances, this means that the optimized bike network is expected to remain close-to-optimal with respect to the all-pairs objective, even under changing demand patterns.

\FloatBarrier

\vspace{1em} \subsection{Effectiveness of spatial edge set reduction}\label{sec:synth_kpath}

In a similar spirit, we analyze the increase of the objective value when reducing the considered edge set to a spatial corridor around the shortest path. As explained in \autoref{sec:method_kpaths}, the runtime can be reduced by allowing flow only on a subset of the edges per s-t-path, denoted $E_H(s,t)$. By measuring the objective value for solutions generated with and without this relaxation, we aim to quantify the effect of the relaxation on the solution quality. We generate twelve graphs of size $n\in \{50, 100, 150, 200\}$ (three per size), and first solve the problem with the proposed simplification, considering only $E_H(s,t)$ instead of $E'$ as flow variables per s-t-path. We thereby set the parameter $\eta$ to 20, 30, 40 and 50 (see \autoref{sec:method_kpaths}). Each solution of the LP yields capacities $\lambda^*(G, \eta)$ that are optimal with respect to the input. For small $\eta$, we expect the solution to diverge substantially from the optimal solution, whereas for higher $\eta$, we hope to obtain similarly good solutions. For a fair comparison, we solve the problem again fixing the capacities to $\lambda^*(G, \eta)$, thereby computing the minimal objective value under the capacity constraints. This objective value is compared to the optimal objective value of the unrestricted problem, where all edges are considered.

\begin{figure}
    \centering
    \begin{subfigure}[b]{0.48\textwidth}
    \includegraphics[width=\textwidth]{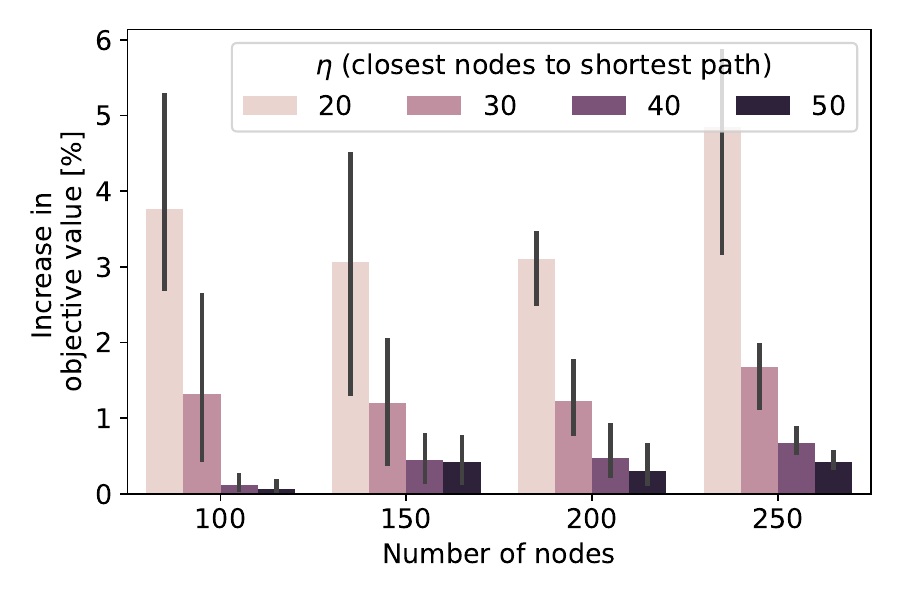}
    \caption{Objective increase}
    \label{fig:kshortest_objective}
    \end{subfigure}
    \hfill
    \begin{subfigure}[b]{0.48\textwidth}
    \includegraphics[width=\textwidth]{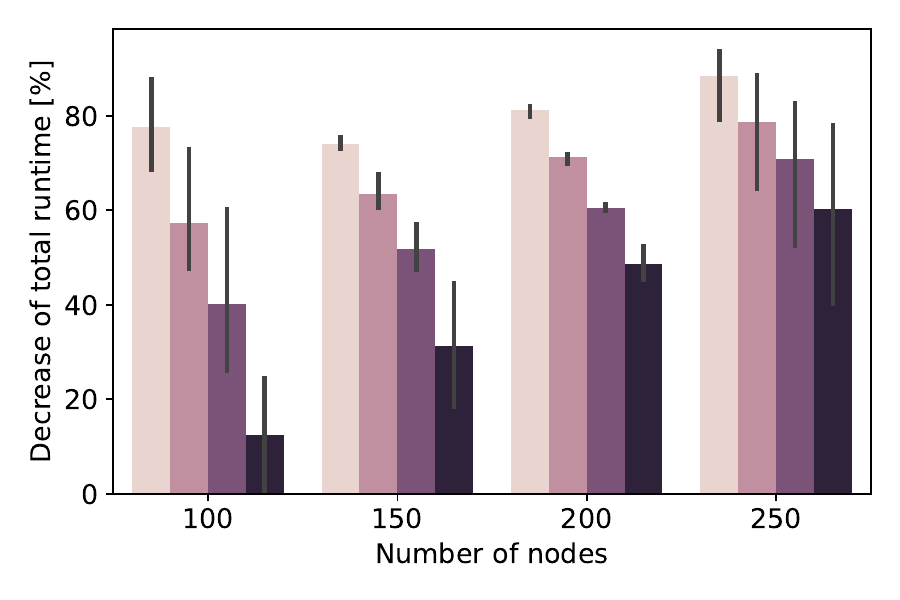}
    \caption{Runtime decrease}
    \label{fig:kshortest_runtime}
    \end{subfigure}
    \caption{Evaluating runtime benefits and reduction in solution quality for restricting the evaluated edges to the spatial neighborhood around the shortest path. While the runtime decreases significantly for larger graphs, the objective value increases by less than 1\% for sufficiently large $k$}
    \label{fig:kshortest_comp}
\end{figure}

The results are shown in \autoref{fig:kshortest_comp}. Indeed, considering less edges significantly improves the runtime, but comes at the cost of slightly worse solutions. However, setting $\eta=40$ or $\eta=50$ yields a good trade-off between efficiency and effectiveness, since the runtime decreases by up to 60\% whereas the objective value increases by less than 1\%. Most importantly, the effect on the runtime becomes stronger with the graph size, whereas the difference in objective remains rather constant, proving the value of this method for reducing the runtime for large real-world instances.

\FloatBarrier

\section{Discussion and conclusion}\label{sec:discussion}

Bike lane allocation is gaining importance as many cities reconsider the priority for cars in the status quo and move towards redesigning cityscapes for active mobility. Our mathematical formalization of this task provides a framework for developing scalable and effective methodologies. % There is a compelling need for methods that provide globally optimal solutions and minimize the impact on the car accessibility and public transport. 
% summary of contributions
Most importantly, our study introduces a novel perspective on bike network planning by prioritizing the trade-off between car and bike travel times through Pareto optimality.
In contrast to most previous work on bike lane planning, we call for an increased emphasis on the impact of bike lanes on the car network, specifically to quantify the trade-off between bike and car travel times. This aim is formalized in the bike lane allocation problem (BNAP) and addressed with our IP formulation that yields Pareto-optimal solutions. To enable the application to real-world scenarios, we developed several relaxations and post-processing schemes, while also integrating key innovations such as demand-driven aspects and the assessment of perceived bike travel times. In contrast to purely heuristic approaches, our solutions are based on the output of a LP, thus maintaining a global perspective of the network and avoiding bottlenecks. At the same time, our framework incorporates advantages of previous bike lane planning methods, such as considering travel demand data and yielding a sequence of networks instead of a single solution.

% results of experiments
The experiments on synthetic data demonstrate the algorithm's efficiency, delivering solutions for networks with hundreds of edges within minutes. 
Our approach also benefits from certain problem relaxations, which significantly accelerate processing times with only slight impacts on travel times. 
The results on real-world instances support the superiority of our approach over heuristic methods. These iterative methods based on the betweenness centrality are inferior in terms of  Pareto optimality due to the lack of integrated perspective of car and bike travel.
Moreover, a detailed case study in Zurich showcased the algorithm's scalability and effectiveness in managing extensive street networks through strategic partitioning and parallelization, further underscoring the utility our method brings to the field of urban planning.

% limitations: shared lane penalty factor
There are several opportunities to enhance the presented approach. For instance, our experiments employ a simplistic model of the perceived travel time of cyclists, applying a uniform penalty for cycling on a car lane. As suggested by \citet{steinacker_demand-driven_2022}, this factor could vary based on the type of adjacent car lane, although our model assumes that bike lanes are properly separated from car lanes. Enhancing network connectivity represents another research direction. Bike network connectivity could be incorporated in our problem formulation by penalizing discrepancies between bike flow $f^b$ and shared flow $f^{\beta}$ at junctions. 

% limitations: traffic flow
Addressing actual traffic flow remains a significant hurdle. Accurately modeling traffic flow typically requires sophisticated simulators, making it challenging to employ optimization algorithms without resorting to a bi-level formulation. 
However, ignoring traffic flow has certain effects on the optimal solution. For example, our algorithm is predisposed to allocate one lane of any double-lane road to bicycles, not accounting for the impact on car travel times due to reduced road width. 
While our study introduces a method to factor in traffic-dependent flow through $\phi$ (see \autoref{eq:phi_traffic_flow}), it remains unclear how to scale the capacity constraints. By setting $\phi=1$ in our experiments and bounding the capacity by the number of lanes, the algorithm is guaranteed to yield feasible solutions, but does not account for traffic congestion~\citep{kumar2021applications} and general flow dynamics~\citep{kessels2019traffic}.% at the cost of limited ability to model traffic jams and traffic flow in general. 

Future work should thus focus on accurately modeling traffic flow and on  adapting our algorithm to account for factors critical to real-world traffic dynamics, such as intersection layouts. Additionally, there are straightforward extensions, such as incorporating parking space allocation constraints or incentivizing bike lanes near green spaces, that could further refine the model. 
To make our algorithm more accessible, developing a user interface that allows parameter adjustments through sliders and visualizes the resultant street networks is crucial. Such a tool would empower urban planners with limited technical expertise to leverage our algorithm for informed decision-making.
% \vspace{2em} \section{Conclusion}\label{sec:conclusion}
%Our algorithm has demonstrated superiority over heuristic methods, proving its applicability to real-world scenarios. Furthermore, its flexible design allows for the integration of additional real-world considerations such as the reallocation of parking space or the proximity of bike lanes to green spaces.
Through this work, we aim to catalyze further research that builds upon our multi-modal viewpoint, aspiring towards a holistic optimization framework for comprehensive bike network planning.

\section*{Acknowledgments}
We thank Aurélien Borgeaud for his contribution in the context of a semester project. This project is part of the E-Bike City project funded by the Department of Civil and Environmental Engineering (D-BAUG) at ETH Zurich and the Swiss Federal Office of Energy (BFE).
The second author was funded by the European Research Council (ERC) under the European Union's Horizon 2020 research and innovation programme (grant agreement No 817750).

\begin{flushright}
    \includegraphics[height=13mm]{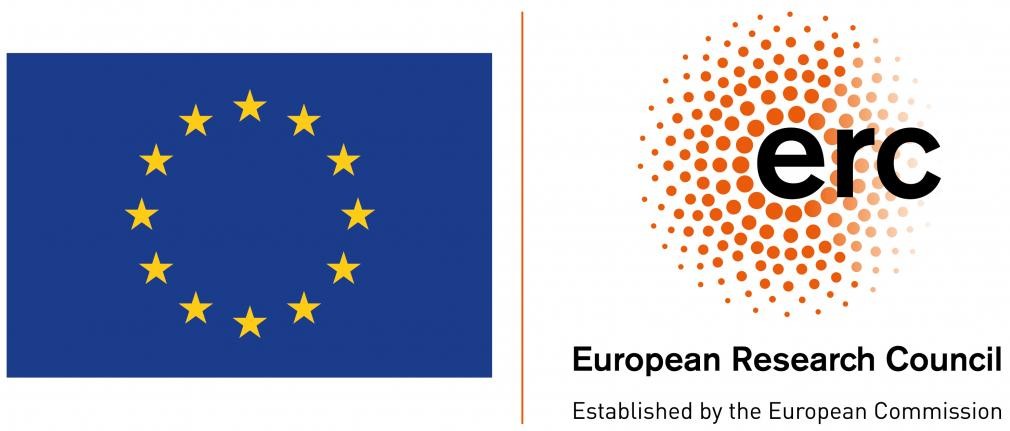}
\end{flushright}

%Bibliography
\section*{Author contributions}
% check CRediT (Contributor Roles Taxonomy) : https://credit.niso.org/ 
Nina Wiedemann: Conceptualization, Data curation, Formal Analysis, Methodology, Software, Validation, Visualization, Writing – original draft\\
Christian Nöbel: Conceptualization, Formal Analysis, Methodology, Software, Writing – review \& editing\\
Lukas Ballo: Data curation, Investigation, Writing – review \& editing\\
Henry Martin: Conceptualisation, Visualization, Writing – review \& editing\\
Martin Raubal: Funding acquisition, Project administration, Supervision, Writing – review \& editing

\bibliographystyle{plainnat}
\bibliography{references}

% \newpage
\vspace{5em}
\appendix
\section*{Appendix}

% \section{Hypervolume indicator}\label{app:hi}

% \begin{figure}[htb]
%     \centering
%     \includegraphics[width=0.5\textwidth]{figures/hi_explained.pdf}
%     \caption{Explanation of computing the inverse HI}
%     \label{fig:hi_explained}
% \end{figure}

\section{Preprocessing of real network data}\label{sec:preprocessing}

The SNMan library builds up on the \texttt{networkx} and \texttt{OSMnx}~\citep{boeing2017osmnx} packages and first constructs a \textit{street} graph with one node per intersection and one edge per street. It further provides functionality to convert the street graph into the \textit{lane graph}, a directed multigraph. Each node in the lane graph is defined by geographic coordinates and elevation, and the edges are enriched with attributes for their distance, speed limit and the type of lane based on available OSM data, allowing to derive $\delta(e)$, $t^c(e)$ and $t^b(e)$. 

An origin-destination $\Omega$ expressing real-world travel demand is derived from public bike sharing or census data. For districts in Zurich, we take all trips in the Mobility Microcensus from Switzerland that intersect with the district region. The Mobility and Transport Microcensus is a statistical survey on travel behavior that is published by the Federal Office for Spatial Development. Overall, it contains travel survey data for more then 57k participants. After intersecting the origin-destination-lines with the city district, 1061 and 508 trips remain for Birchplatz and Affoltern respectively. The trip origin and destination are matched to nodes in the graph by selecting the node closest to their geographic coordinates, yielding 498 and 267 unique OD-pairs respectively (see \autoref{tab:instances}).

On the other hand, to the best of our knowledge there is no public travel data for Chicago or Cambridge. Instead, we utilize public data from their respective local bike sharing services. It is worth noting that this biases the OD-data to cycling movement; however, it can be assumed that the data still reflect typical mobility behavior within the district, irrespective of the transport mode. For Chicago, data is available from Divvy bike sharing\footnote{\url{https://divvybikes.com/system-data }}, whereas in Cambridge, the operating service is Bluebikes\footnote{\url{https://bluebikes.com/system-data}}. In both cases, we download and merge all trip data from the whole year of 2023 and only filter out stations that are marked as charging or maintenance stations. The origin-destination data is intersected with the respective district based on their beeline connection, and the resulting subset of OD-pairs is mapped to graph nodes by their geographic locations, yielding a set of node-based OD-pairs. Due to the large number of resulting pairs, we select only the OD-pairs that collectively account for 75\% of the trips from 2023, resulting in the counts listed in \autoref{tab:instances}.

\section{City partitioning for case study}\label{app:partitioning}

\begin{figure}[bht]
    \centering
    \begin{subfigure}[b]{0.49\textwidth}  
    \includegraphics[width=\textwidth]{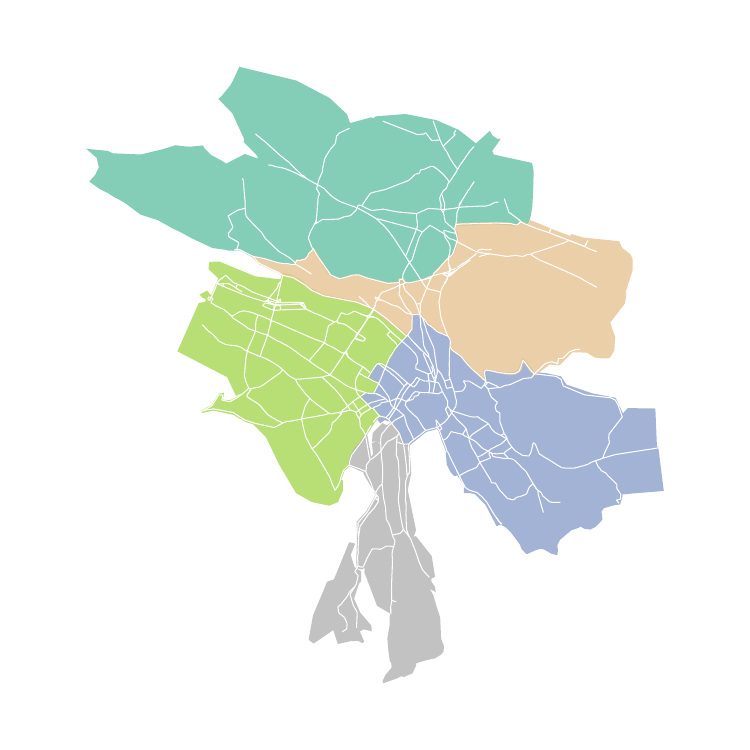}
    \caption{Main road division}
    \label{fig:main_partition}
    \end{subfigure}
    \begin{subfigure}[b]{0.49\textwidth}  
    \includegraphics[width=\textwidth]{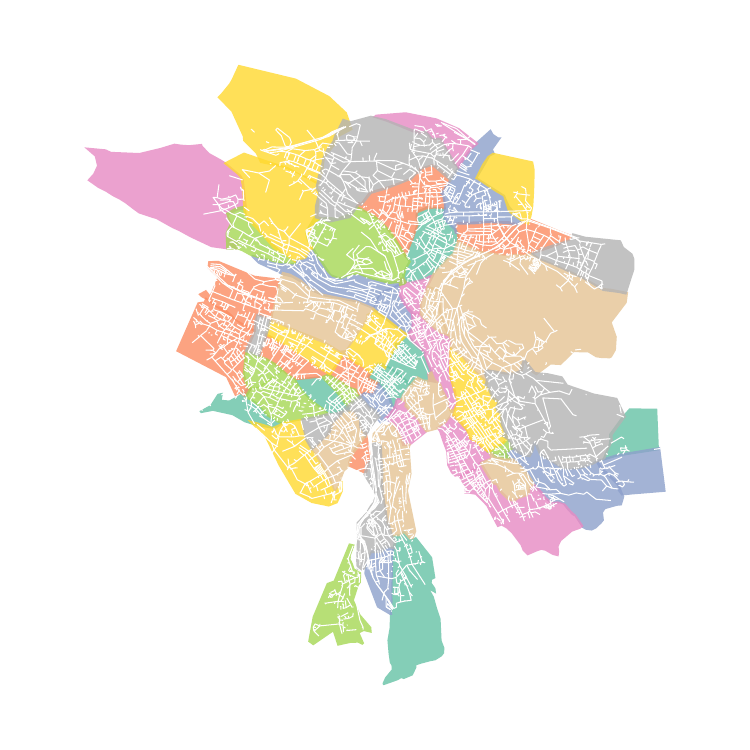}
    \caption{Non-main-road division}
    \label{fig:side_partition}
    \end{subfigure}
    \caption{Partitioning of Zurich to rebuild the full city street network. First, the main roads are partitioned and bike lanes are allocated. Secondly, the side roads are partitioned by districts and bike lanes are allocated.}
    \label{fig:partitioning}
\end{figure}

We devide the city of Zurich into 57 parts following \citet{ballo2023}, as shown in \autoref{fig:side_partition}. Since the edges are still defined between \textit{intersections}, the main-road graph still comprises 1902 edges. The main roads are thus partitioned as well, by a manual aggregation of postal codes (see \autoref{fig:main_partition}). When optimizing the side roads, main roads are not filtered out but are fixed. Highways as well as lanes allocated for public transport are not changed in both cases.

\FloatBarrier

\section{Baselines}\label{sec:baselines_app}

\begin{figure}[htb]
    \centering
    \includegraphics[width=\linewidth]{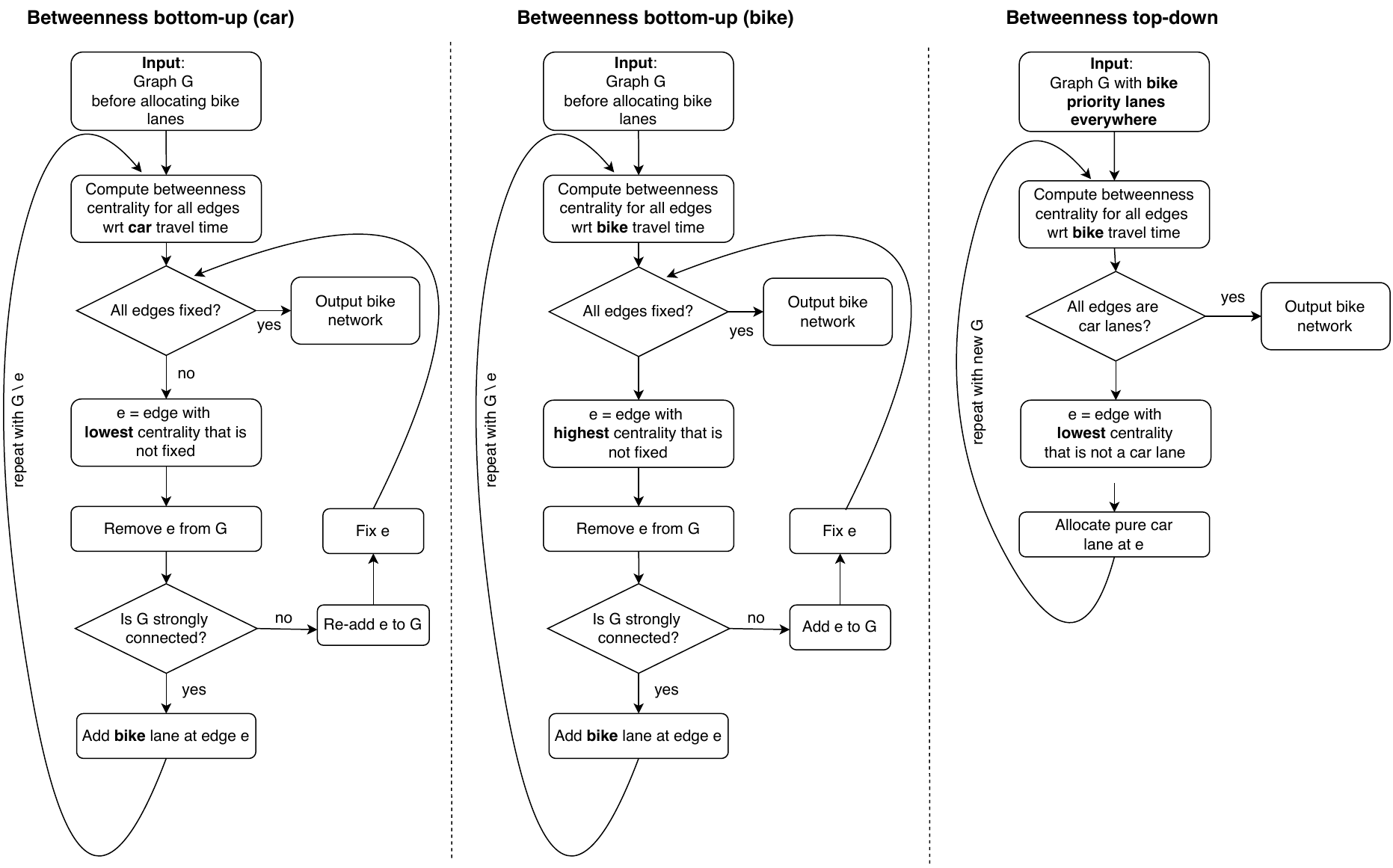}
    \caption{Flowcharts of the three baseline algorithms. All algorithms are heuristic algorithms based on the betweenness centrality. However, they take different strategies wrt. the prioritization of bikes and cars (see bold differences).}
    \label{fig:baseline_flowchart}
\end{figure}

\autoref{fig:baseline_flowchart} provides details on the three baseline algorithms. They use betweenness centrality as a heuristic and iteratively re-allocate edges to bikes or cars, while testing for strong connectivity of the car network in every step. They differ in the prioritization of bikes and cars, either allocating bike lanes at the most important edges for cyclists, or at the least important edges for car drivers.

\FloatBarrier
\section{Runtime of Integer Program}\label{sec:integer_runtime}

\begin{figure}[htb]
    \centering
    \includegraphics[width=0.5\linewidth]{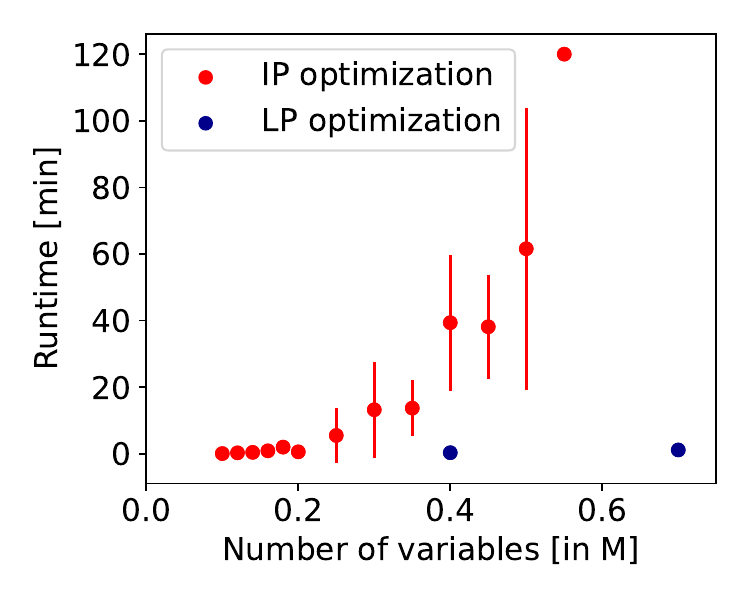}
    \caption{Runtime comparison between IP and LP. The IP is order of magnitudes slower than the LP on instances of more than 200k nodes.}
    \label{fig:integer_runtime}
\end{figure}

\autoref{fig:integer_runtime} shows how the runtime of the IP scales with the instance size. For instances of 500k variables, where the LP can still be solved within minutes, solving the IP already takes one hour on average. Considering the size of real-world instances, the LP relaxation is usually necessary.

\section{Comparison of LP and IP solution}\label{app:int_vs_lin}

\begin{figure}
    \centering
    \begin{subfigure}[b]{0.48\textwidth}
    \includegraphics[width=\linewidth]{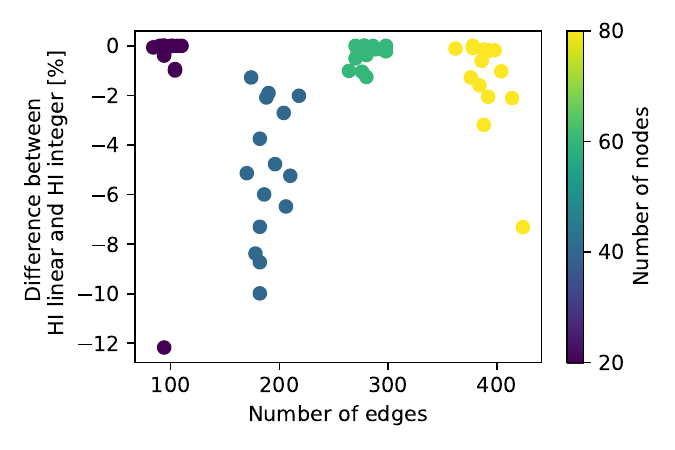}
    \end{subfigure}
    \hfill
    \begin{subfigure}[b]{0.48\textwidth}
    \includegraphics[width=\linewidth]{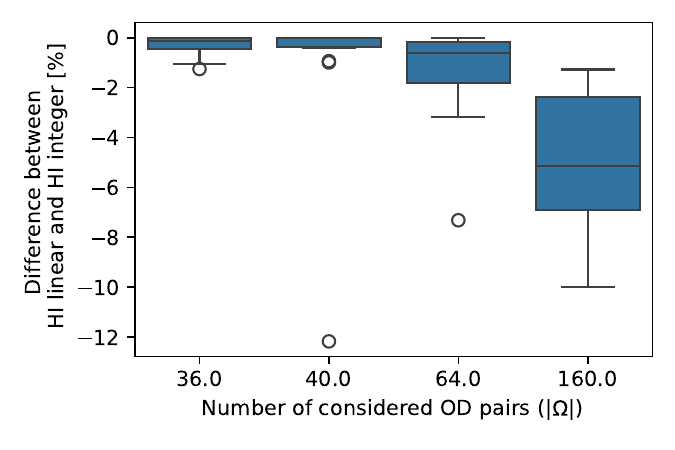}
    \end{subfigure}
    \caption{Dependence of the relative difference between LP and IP solution on the graph size and OD size. While there is no clear impact of the number of edges or graph density, the number of considered OD-pairs plays a crucial role.}
    \label{fig:int_vs_lin_od}
\end{figure}

The results presented in \autoref{sec:synthetic_lp_vs_ip} showed that the LP and rounding yields solutions that are close to the IP solution in most cases. To understand the conditions for achieving the optimal solution, we analyze the effect of certain network properties  on the optimality in \autoref{fig:int_vs_lin_od} and \autoref{tab:regression_results}. \autoref{fig:int_vs_lin_od} shows that the divergence between linear and integer HI increases with the number of considered OD pairs, whereas the graph size itself does not have a strong effect. \autoref{tab:regression_results} provides a statistical analysis of the effect size of network properties in a linear regression model (normalized for comparability), where the dependent variable is the relative difference between LP and IP HI. The effect of the number of OD pairs is significant. In addition, the standard deviation of the clustering coefficient of the graph seems to be negatively associated with optimality. Intuitively, graphs with high structural variances within themselves are more likely to have bottlenecks where the relaxation fails. The positive relation between the mean clustering coefficient and the optimality is likely due to clustered areas offering multiple viable options for bike lane placement.

\begin{table}
\centering
\caption{Regression Coefficients for the effect of network properties on the difference between LP and IP solution}
\label{tab:regression_results}
\begin{tabular}{lrr}
\toprule
 & Coefficient & p-value \\
\midrule
Intercept & -1.9292 & 0.0000 \\
Number of nodes & 2.3454 & 0.3638 \\
Number of edges & -2.9632 & 0.2741 \\
Number of OD-pairs ($|\Omega|$) & -1.9892 & 0.0000 \\
Mean betweenness centrality & -2.2414 & 0.2008 \\
Standard deviation of betweenness centrality & 1.3406 & 0.2208 \\
Mean clustering coefficient & 1.3995 & 0.0845 \\
Standard deviation of clustering coefficient & -1.4919 & 0.0446 \\
\bottomrule
\end{tabular}
\end{table}

\end{document}